\documentclass[10pt,journal,compsoc]{IEEEtran}

\ifCLASSOPTIONcompsoc
  \usepackage[nocompress]{cite}
\else
  \usepackage{cite}
\fi

\usepackage{amssymb}
\usepackage{amsfonts}
\usepackage{amsmath}
\usepackage{graphicx}%
\usepackage{verbatim}

\usepackage{algorithm}
\usepackage{algpseudocode}
\usepackage{subeqnarray}
\usepackage{amsmath,bm}
\usepackage{upgreek}
\usepackage{hyperref}
\usepackage{amsthm}
\usepackage{setspace}
\usepackage{subeqnarray}
\usepackage[mathscr]{eucal}
\usepackage{subfig}
\usepackage{bm}
\usepackage{lipsum}
\makeatletter
\def\footnoterule{\relax%
	\kern-3pt
	\hbox to \columnwidth{\hfill\vrule width 0.5\columnwidth height 0.4pt\hfill}
	\kern4.6pt}
\makeatother

\allowdisplaybreaks[4]
\newcommand*{\rom}[1]{\expandafter\@slowromancap\romannumeral #1@}
\newtheorem{theorem}{Theorem}

\newtheorem{definition}{Definition}
\newtheorem{lemma}{Lemma}

\newtheorem{remark}{Remark}

\algnewcommand{\Inputs}[1]{%
	\State \textbf{Inputs:}
	\Statex \hspace*{\algorithmicindent}\parbox[t]{.9\linewidth}{\raggedright #1}
}
\algnewcommand{\Initialize}[1]{%
	\State \textbf{Initialization:}
	\Statex \hspace*{\algorithmicindent}\parbox[t]{.9\linewidth}{\raggedright #1}
}

\begin{document}
\title{Optimization of convergence rate via algebraic connectivity}

\author{Zhidong He~\IEEEmembership{}
\IEEEcompsocitemizethanks{\IEEEcompsocthanksitem Zhidong He is with the Faculty of Electrical Engineering, Mathematics and Computer Science, Delft University of Technology, The Netherlands.
P.O Box 5031, 2600 GA Delft.\protect\\
E-mail: Z.He@tudelft.nl}
\thanks{}}


\IEEEtitleabstractindextext{%
\begin{abstract}
The algebraic connectivity of a network characterizes the lower-bound of the exponential convergence rate of consensus processes. 
This paper investigates the problem of accelerating the convergence of consensus processes by adding links to the network. Based on a perturbation formula of the algebraic connectivity, we propose a greedy strategy for undirected networks and give a lower bound of its performance through an approximation of submodularity.
We further extend our investigation to directed networks, where the second smallest real part among all the eigenvalues of the non-Hermitian Laplacian matrix, i.e., the generalized algebraic connectivity, indicates the expected convergence rate. We propose the metrics to evaluate the impact of an adding subgraph on the generalized algebraic connectivity, and apply a modified greedy strategy to optimize the generalized algebraic connectivity.
Numerical results in empirical networks exhibit that our proposed methods outperform some other methods based on traditional topological metrics. Our research also verifies the dramatic differences between the optimization in undirected networks and that in directed networks.
\end{abstract}

\begin{IEEEkeywords}
convergence rate, algebraic connectivity, consensus process, directed network
\end{IEEEkeywords}}
\maketitle

\IEEEdisplaynontitleabstractindextext
\IEEEpeerreviewmaketitle

\IEEEraisesectionheading{\section{Introduction}\label{sec:introduction}}

%
%
%
%
\IEEEPARstart{A} large number of collective dynamic processes, e.g., the Markovian process \cite{van2014performance}, the network diffusion \cite{simonsen2004diffusion}, the fluid flow in tank systems \cite{luyben1989process}, the electric network \cite{van2017pseudoinverse}, has a similar objective that all nodes reach an agreement regarding a certain quantity of interest by exchanging the nodal states with their neighboring nodes. These dynamics can be generally described by the consensus model in networks \cite{olfati2007consensus,li2010consensus}.
The convergence rate of a consensus process on a network indicates the speed that each node tends to its final steady state (assuming that the steady state exists), which characterizes the behavior of an autonomous system and can be an important metric of the performance (e.g. robustness or efficiency) for networked control systems.

The behavior of the consensus model can be featured by the eigensystem of the Laplacian matrix of the adjacency matrix, where the algebraic connectivity of the undirected network indicates a lower-bound of the exponential convergence rate \cite{olfati2007consensus}. Thus, the problem of convergence acceleration of a process reduces to increasing the algebraic connectivity of the undirected network by topological adjustments. If the variables of interconnections are continues-valued, i.e., the links are weighted, maximizing the algebraic connectivity constrained by a total budget of link weights reduces to a semi-definite programming \cite{ghosh2006growing}, which can be solved by a subgradient method \cite{kim2006maximizing,boyd2004fastest}. 
Meanwhile, several heuristic methods based on the topological metrics are proposed to approach the NP-hard integer-optimization problem in the case of unweighted networks \cite{sydney2013optimizing,ogiwara2017maximizing,alenazi2013network,li2018maximizing}, but the performance is not guaranteed.

Moreover, several flow processes, e.g., the impedance circuit, the liquid flow, the traffic flow, yields the fact that the interactions between nodes are directed instead of bidirected. The convergence rate to consensus in a directed network depends on both the topology of the network and its initial state vector \cite{li2010cooperative}. The expected convergence rate (ECR) is defined to measure the speed of convergence for random initial state vectors in directed networks \cite{asadi2017expected}.
Notwithstanding the importance and generality of the expected convergence rate in directed networks, the method of accelerating consensus processes in directed networks by topological adjustments, to the best of our knowledge, has seldom been investigated. 

In this paper, we investigate the impact of the network perturbation on the convergence rate of dynamics via the algebraic connectivity both in undirected and directed networks, which are confined to unweighted networks.
For undirected networks, we present a perturbation formula of the algebraic connectivity for an adding subgraph, which shows the impact of multiple adding links on the algebraic connectivity, as well as demonstrates that the algebraic connectivity resembles an approximately submodular function with respect to the adding subgraphs. Based on the metric for the impact of an individual link, the greedy strategy to increase the algebraic connectivity by adding links one by one could have a constant performance.

For directed networks, we confine ourselves to the strongly connected network without self-loops to ensure the existence of the non-trivial steady state \cite{li2010consensus}. We show that the second smallest real part among all the eigenvalues of the Laplacian matrix indicates the lower bound of the expected convergence rate in directed networks, which is defined as the generalized algebraic connectivity \cite{asadi2017expected}. 
Based on the bounds of the generalized algebraic connectivity perturbation with respect to the adding subgraph, we propose the metrics to measure the impact of an individual link and heuristically apply the greedy strategy to maximize the generalized algebraic connectivity.
Further, we compare the performance of the proposed methods with other heuristic methods based on traditional nodal centralities, and discuss the difference between the strategies for increasing the convergence rate between undirected and directed networks.

The main contribution of this paper can be summarized as:
\begin{itemize}
	\item[1.] We apply a greedy strategy to maximize the algebraic connectivity $\mu$ in undirected networks, which guarantees a constant performance.
	\item[2.] Based on the bounds of the generalized algebraic connectivity $\Re(\mu)$ perturbation, we propose heuristic greedy strategies to maximize $\Re(\mu)$ for directed networks.
	\item[3.] Numerical tests in real-world networks show the superiority of our proposed methods compared with some other heuristics.
\end{itemize}

The remainder of this paper is organized as follows. The physical significance of the algebraic connectivity and the definition of the generalized algebraic connectivity are introduced in Section 2. Section 3 provides a perturbation formula of the algebraic connectivity and propose the greedy algorithm in undirected networks. We extend the optimization problem and propose some heuristic strategies for directed networks in Section 4. We evaluate the proposed methods in Section 5 and conclude this paper in Section 6.

\section{Algebraic connectivity for consensus processes in networks}
\subsection{Consensus processes in undirected networks}
We consider that an undirected and unweighted network $G(\mathcal{N},\mathcal{L})$ consisting of set $\mathcal{N}$ with $N$ nodes and set $\mathcal{L}$ with $L$ links is represented by the adjacency matrix $A$. The entry of the adjacency matrix $a_{ij} = 1$ if there is a link between node $i$ and node $j$, and otherwise $a_{ij} = 0$.
Let $v_i(t)$ denotes the state of node $i$ at time $t$, 
the consensus model follows $\frac{d v_i(t)}{dt} = \sum_{j=1}^N a_{ij}(v_i(t)-v_j(t))$, and the vector form for the state vector $\bm{v}(t) = (v_1(t),v_2(t),\dots, v_N(t))^T$ follows
\begin{align}
\frac{d \bm{v}}{dt} = -(\Delta-A)\bm{v}(t) = -Q\bm{v}(t)
\end{align}
where $\Delta = \text{diag}(d_1,d_2,\dots,d_N)$ is the diagonal degree matrix (e.g., $d_i$ denotes the degree of node $i$), and $Q$ is the Laplacian matrix of $A$.

We assuming that all the eigenvalues $\lambda_k$ of $Q$ are distinct, i.e., $\lambda_1<\lambda_2<\dots<\lambda_N$, associated with the eigenvector $x_1,x_2,\dots,x_N$.
The state vector $\bm{v}(t)$ can be written as $\bm{v}(t) = e^{-Qt}\bm{v}(0)$ with the initial state vector $\bm{v}(0)$.
The smallest eigenvalue in the symmetric matrix $Q$ is equal to $\lambda_N(Q) = 0$, which corresponds to the steady-state eigenvector $\bm{\pi} = x_N =\bm{u}^T$.
Thus, the state vector can be written by 
\begin{align}
\bm{v}(t)=\bm{\pi}+\sum_{k=2}^N e^{-\lambda_k t} x_k x_k^T \bm{v}(0)
\end{align} 
which implies that the state $v_i(t)$ of each node exponentially converges to the steady state. The relation 
\begin{align}
-\lim_{t\rightarrow\infty}\frac{\log |v_i(t)-\pi_i|}{t} \geq \min_{k\leq N-1}\{\lambda_k\}
\end{align}
implies that the second smallest eigenvalue $\lambda_{N-1}$ of the Laplacian matrix $Q$, i.e., the algebraic connectivity, indicates a lower-bound of the expected exponential convergence rate for each nodal state with a random initial state.

\subsection{Consensus processes in directed networks}
We further consider the dynamic processes in the directed network $G(\mathcal{N},\mathcal{L})$ with $N$ nodes and $L$ directed links. The entry $a_{ij} = 1$ in the asymmetric adjacency matrix $A$ represents a directed link with source node $i$ and target node $j$.
In the consensus process, the existence of link $\ell_{ji}$ represents the communication from node $i$ to node $j$. Thus, the consensus model follows $\frac{d v_i(t)}{dt} = \sum_{j=1}^N a_{ji}(v_j(t)-v_i(t))$, and the state vector follows
\begin{align}\label{equ:consensus_model_directed}
\frac{d \bm{v}(t)}{dt}=-(\Delta_{in}-A^T)\bm{v}(t)
\end{align}
where $\Delta_{in}$ is the diagonal nodal in-degree matrix. We define the operator $Q=\Delta_{in}-A^T$ as the generalized Laplacian matrix in directed network\footnote{Some researchers \cite{olfati2007consensus} prefer applying link $\ell_{ji}$ to represent communication from node $i$ to node $j$, then the Laplacian matrix is $Q = \Delta_{out}-A$ with the diagonal nodal out-degree matrix $\Delta_{out}$. Our definition in this paper has no difference with the previous for system analysis, but can be more intuitive for topological changing.}.
If the network is strongly connected \cite{olfati2007consensus}, i.e., there exists a path from one node to another for any two nodes, the smallest eigenvalue $\lambda_N$ of $Q$ is equal to 0, and the system can reach a steady state $\bm{\pi}$.
Assuming that the generalized Laplacian matrix has distinct eigenvalues $\lambda_k$ and corresponding right- and left- eigenvectors $x_k$ and $y_k$, the state vector can be generalized as 
\begin{align}
\bm{v}(t)=\bm{\pi}+\sum_{k=2}^N e^{-\lambda_k t} x_k y_k^T \bm{v}(0)
\end{align} 
Since the complex eigenvalue for the non-Hermitian matrix $Q$ can be written by $\lambda_k = \Re(\lambda_k)+i\Im(\lambda_k)$ with non-negative real parts $\Re(\lambda_k)$, we rewrite the state vector by $\bm{v}(t)=\bm{\pi}+\sum_{k=2}^N e^{-\Re(\lambda_k)t-i\Im(\lambda_k)t } x_k y_k^T \bm{v}(0)$. 
Denoting $\bm{e}_i$ the basis vector, the difference between the nodal state $v_i(t)$ and its steady state $\pi_i$ follows
\begin{small}
\begin{align}
|v_i(t)-\pi_i|&=\bm{e}_i\bigg|\sum_{k=2}^N e^{-\Re(\lambda_k)t-i\Im(\lambda_k)t } x_k y_k^T \bm{v}(0)\bigg| \notag\\
&\leq\sum_{k=2}^N e^{-\Re(\lambda_k)t} |\bm{e}_i x_k y_k^T \bm{v}(0)|  
\end{align}
\end{small}
which yields that the expected exponential convergence rate has a lower bound
\begin{align}
-\lim_{t\rightarrow\infty}\frac{\log |x_i(t)-\pi_i|}{t} \geq \min_{k\leq N-1}\{\Re(\lambda_k)\}
\end{align}

For brevity, we write the algebraic connectivity as $\mu = \lambda_{N-1}$ in undirected networks, and define the generalized algebraic connectivity $\Re(\mu) = \min_{i<N}\Re(\lambda_i)$ as the second smallest real part among all the eigenvalues of the Laplacian matrix in directed networks.

\section{Algebraic connectivity in undirected networks}
In this section, we derive an approximation of the algebraic connectivity for the topological perturbation in undirected networks. Further, we propose a greedy method to maximize the algebraic connectivity by adding links or subgraphs.

\subsection{Topological perturbation for the algebraic connectivity}
We denote by $A$ the adjacency matrix of the original undirected network $G(\mathcal{N},\mathcal{L})$, and by $A+\Updelta A$ the adjacent matrix of the network under perturbation $\Updelta A$. The $N\times N$ perturbation matrix $\Updelta A$ is also symmetric, which is positive-definite for the adding subgraph, and negative-definite for deleting a subgraph. 
The $N\times N$ Laplacian matrix $\Updelta Q$ corresponding to $\Updelta A$ is the perturbation to the original Laplacian matrix $Q$, which is symmetric and positive-definite for an adding subgraph (negative-definite for a delated subgraph).
According to the eigenvalue perturbation theorem \cite{van2010graphspectra, stewart1990matrix}, the perturbation of the algebraic connectivity follows
\begin{small}
\begin{align}\label{equ:asmptotic_expan}
\Updelta \mu = x_{N-1}^T\Updelta Qx_{N-1}+\sum_{k\neq N-1}^N\frac{(x_k^T\Updelta Qx_{N-1})^2}{\lambda_{N-1}-\lambda_k}+\mathcal{O}(\Vert\Updelta Q\Vert_F^3)
\end{align}
\end{small}
where $x_{N-1}$ is the Fiedler vector and $\Vert\cdot\Vert_F$ is the Frobenius norm. However, estimating the algebraic connectivity perturbation $\Updelta\mu$ by \eqref{equ:asmptotic_expan} requires all the eigenvalues and the eigenvectors, whose computational complexity is usually high for large networks, i.e. $\mathcal{O}(N^3)$ by QR decomposition.

We hereby present an approximate of the algebraic connectivity for the topological perturbation $\Updelta A$. The Laplacian matrix of a network can be convert into a row-stochastic and non-negative matrix $S$, by using the transformation $S=I_N - \epsilon Q$, where $I_N$ is the identity matrix and $\epsilon>0$ is a sufficiently small number. The matrix $S$ is irreducible if the network $G$ is connected when the parameter $\epsilon < \frac{1}{d_{\max}}$ ensures $S$ is non-negative \cite{olfati2007consensus}.
We denote by $\gamma_1$ the left eigenvector of $S$ corresponding to its largest eigenvalue 1.
Then, we can construct the matrix $R = S-\frac{\gamma_1\gamma_1^T}{\Vert\gamma_1\Vert^2}$, and denote by $z$ the left eigenvector of $R$ corresponding to the largest eigenvalue $\lambda_1(R)$. The algebraic connectivity $\mu$ of $Q$ shifts to the largest eigenvalue of $R$, and can be computed \cite{li2013distributed, quarteroni2008matematica}  by
\begin{align}
\mu(Q)  = \frac{1}{\epsilon}(1-\lambda_2(S)) =\frac{1}{\epsilon}(1-\lambda_1(R))
\end{align}

We suppose that the algebraic connectivity $\mu$ becomes $\widetilde{\mu}=\mu+\Updelta\mu$ if the Laplacian matrix of a network $Q$ becomes $Q+\Updelta Q$.
The matrix $S$ after perturbation is $S+\Updelta S$ with $\Updelta S = -\epsilon\Updelta Q$.
The eigenvector $\gamma_1$ corresponding to $\lambda_1(S)$ normalized by $\Vert\gamma_1\Vert^2 = 1$ is equal to the eigenvector of $Q$ corresponding to $\lambda_N(Q) = 0$, i.e., $\gamma_1(S)=x_N(Q) = (\frac{1}{\sqrt{N}},\frac{1}{\sqrt{N}},\dots, \frac{1}{\sqrt{N}})^T$. Thus, the eigenvector $\gamma_1$ does not change for the perturbation, i.e., $\gamma_1(S+\Updelta S) = \gamma_1(S)$, which yields the perturbation $\Updelta R = \Updelta S$.

Denoting $z$ the principle eigenvector of the matrix 
\begin{align}\label{equ:defineR}
R = I_N - \epsilon Q-\frac{1}{N}J
\end{align}
where $J$ is an unit matrix, we further estimate the eigenvector $z+\Updelta z$ for a small perturbation $\Vert\Updelta R\Vert_2\ll\Vert R\Vert_2$ by means of one further iteration of the power method \cite{milanese2010approximating} as
\begin{small}
\begin{align} 
z+\Updelta z &\approx \frac{(R+\Updelta R)z}{\Vert(R+\Updelta R)z\Vert_2} =
\frac{\lambda_1(R)z+\Updelta Rz}{\Vert \lambda_1(R)z+\Updelta Rz \Vert_2} \notag\\
&= \frac{\lambda_1(R)z+\Updelta Rz}{\sqrt{\lambda_1^2(R)+2\lambda_1(R)z^T\Updelta Rz +z^T\Updelta R^T\Updelta Rz}} \notag\\
&\approx z+\frac{1}{\lambda_1(R)}\Updelta R z = z-\frac{\epsilon}{\lambda_1(R)} \Updelta Q z
\end{align}
\end{small}

We apply a second-order perturbation result as $\Updelta\lambda_1(R)\approx \Updelta \lambda_1^{(1)}(R)+\Updelta \lambda_1^{(2)}(R)$ in discrete calculus.
The first term $\lambda_1^{(1)}(R)$ can be obtained by the first-order eigenvalue perturbation \cite{stewart1990matrix}, i.e., $\Updelta\lambda_1^{(1)}(R) = z^T\Updelta Rz$.
The second term $\Updelta \lambda_1^{(2)}(R)$ following the Taylor theorem equals
\begin{small}
\begin{align*}
\Updelta \lambda_1^{(2)}(R) = \frac{1}{2}\Updelta^2 \lambda_1^{(1)}(R) = \frac{1}{2}\left(\Updelta\lambda_1^{(1)}(R+\Updelta R)-\Updelta\lambda_1^{(1)}(R)\right)
\end{align*}
\end{small}
where the term $\Updelta\lambda_1^{(1)}(R+\Updelta R)$ can be approximated by
\begin{small}
\begin{align}
\Updelta\lambda_1^{(1)}(R+\Updelta R) &= \frac{(z+\Updelta z)^T\Updelta R(z+\Updelta z)}{(z+\Updelta z)^T(z+\Updelta z)} \notag\\
&\approx z^T	\Updelta Rz + z^T\Updelta R\Updelta z +\Updelta z^T\Updelta R z \notag\\
&=z^T(-\epsilon\Updelta Q)z + z^T(-\epsilon\Updelta Q)(-\frac{\epsilon}{\lambda_1(R)} \Updelta Q z) \notag\\
&\quad +(-\frac{\epsilon}{\lambda_1(R)} \Updelta Q z)^T(-\epsilon\Updelta Q)z \notag\\
&= -\epsilon z^T \Updelta Q z + \frac{2}{\lambda_1(R)}z^T(\epsilon\Updelta Q)^2 z
\end{align}
\end{small}

Therefore, we can obtain that the exact perturbation $\Updelta\mu$ of the algebraic connectivity approximates
\begin{small}
\begin{align}\label{equ:2order_perturbation}
\Updelta\mu = -\frac{1}{\epsilon}\Updelta\lambda_1(R) \approx  z^T \Updelta Q z - \frac{\epsilon}{\lambda_1(R)}z^T(\Updelta Q)^2 z
\end{align}
\end{small}
and the algebraic connectivity $\widetilde{\mu}(\Updelta Q)$ after a perturbation $\Updelta Q$ is 
\begin{small}
	\begin{align}\label{equ:mu_perturbation}
	\widetilde{\mu}(\Updelta Q) \approx  \mu+z^T \Updelta Q z - \frac{\epsilon}{\lambda_1(R)}z^T(\Updelta Q)^2 z
	\end{align}
\end{small}

According to the perturbation formula \eqref{equ:mu_perturbation}, we can apply the power iteration method to compute the largest eigenvalue $\lambda_1(R)$ and the principle eigenvector $z$, and further estimate the algebraic connectivity $\widetilde{\mu}$ for the perturbation $\Updelta Q$, which reduces the computational complexity to $\mathcal{O}(N^2)$.
Compared with the first-order perturbation $\Updelta\mu \approx z^T \Updelta Q z =\sum_{i,j\in N} \Updelta a_{ij}(z_i-z_j)^2$ proposed in \cite{ghosh2006growing}, the second term in the derived approximation \eqref{equ:2order_perturbation} exhibits an additional penalized term. Specifically, the second term 
\begin{small}
\begin{align}\label{equ:second_term_str}
z^T(\Updelta Q)^2 z = &2\sum_{i,j\in N} \Updelta a_{ij}(z_i-z_j)^2 \notag \\ &
+ 2\sum_{i,j,k\in N}\Updelta a_{ij}\Updelta a_{jk}(z_i-z_j)(z_k-z_j)
\end{align}
\end{small}
implies the effect of multiple links incident to a same node in a perturbation subgraph.

\subsection{Maximize the algebraic connectivity by adding links}
We then consider the problem of maximizing the algebraic connectivity by adding subgraphs or links.
Ghosh and Boyd \cite{ghosh2006growing} show that maximizing the algebraic connectivity $\mu$ in undirected networks by allocating the continues-valued links weights $w_\ell\in[0,1]$ is a convex problem.
However, the combinatorial optimization problem that maximizing the algebraic connectivity by adding $K$ unweighted links, i.e., $w_\ell\in\{0,1\}$ is a NP-hard problem, whose computational complexity is $\mathcal{O}\left({{\frac{1}{2}N(N-1)-L}\choose{K}}N^3\right)$ by a brute force.
An intuitive and straight way is to apply the greedy strategy to adding links for increasing the algebraic connectivity as Algorithm \ref{alg:org_greedy}, which selects the link that increases the algebraic connectivity most in each iteration.
\begin{small}
	\begin{algorithm}[htb]
		\caption{Original Greedy Method}\label{alg:org_greedy}
		\begin{spacing}{1}
			\begin{algorithmic}[1]
				\Inputs {the undirected network $G$ and the number of links $K$} 
				\Initialize {Set the solution $S$ as an empty link set}
				\While {$|S| \le K$}
				\State {$\ell_{ij} = \arg\max_{\ell\notin G} \{\mu(G\cup \{\ell_{ij}\})\} $}
				\State {$S = S\cup \{\ell_{ij}\}$, $G = G\cup \{\ell_{ij}\}$ }		
				\EndWhile
				\State {Return $S$}
			\end{algorithmic}	
		\end{spacing}	
	\end{algorithm}
\end{small}

\begin{definition} (Submodularity \cite{lovasz1983submodular}): A set function $f$ is submodular if for all subsets $A\subset B\subseteq V$ and all elements $s\notin B$, it holds that
	\begin{align}
	f(A\cup\{s\})-f(A) \geq f(B\cup\{s\})-f(B)
	\end{align}
\end{definition}
Physically, submodularity is a diminishing returns property where adding an element to a smaller set gives a larger gain than adding one to a larger set.	
Greedy strategy is a suggested method of operational simplicity to approach the near-optima for a monotone submodular maximization \cite{feldman2017greed}. 
A celebrated results by Nemhauser \emph{et al.} \cite{nemhauser1978analysis} proves that the greedy method provides a good approximation to the optimal solution of the NP-hard optimization problem. 
Defining the impact of an adding subgraph with the Laplacian matrix $\Updelta Q$ on the algebraic connectivity as 
\begin{align}\label{equ:impact}
	\mathcal{E}(\Updelta Q) = z^T \Updelta Q z - \theta z^T(\Updelta Q)^2 z
\end{align}
where the constant $\theta :=\frac{\epsilon}{\lambda_1(R)}$ depends on the original network, we now investigate the property of the impact function $\mathcal{E}(\Updelta Q)$.

In the following analysis, we assume that the impact function $ \mathcal{E}(\Updelta Q)$ approximates the algebraic connectivity $\Updelta \mu$ well, thus the maximzation of $\mathcal{E}(\Updelta Q)$ and $\Updelta \mu$ are consistent.
Further, computing the exact algebraic connectivity $\mu(G\cup \{\ell_{ij}\})$ in Algorithm \ref{alg:org_greedy} for all candidate links $\ell_{ij}$ in each iteration still requires high computational cost for large scale networks.
According to the impact $\mathcal{E}(\Updelta Q)$ of a subgraph, the impact of an individual link $\ell_{ij}$ can be measured by the metric $\Omega(\ell_{ij}) =|z_i-z_j|$, which can be applied to faster select the best link $\ell_{ij} = \arg\max_{\ell\notin G} \{|z_i-z_j|\}$ for each iteration in Algorithm \ref{alg:greedy_metric}. Coincidentally, the impact \eqref{equ:second_term_str} validates that the second term $z^T(\Updelta Q)^2 z$ suggests the same metric $|z_i-z_j|$ for the importance of an individual link. 
In addition, the improved perturbation function \eqref{equ:impact} with the second term \eqref{equ:second_term_str} can be applied for more general applications, e.g., adding a subgraph with multiple links in each iteration, or batched greedy optimization \cite{liu2018performance}.

The feasibility of the greedy algorithm \ref{alg:greedy_metric} is an open question. Although the algebraic connectivity is a concave function of continuous-valued link weights, Remark \ref{thm:submodular_lambda} shows that the approximated increment of the algebraic connectivity $\mathcal{E}(\Updelta Q)$ by adding unweighted subgraphs is unfortunately non-submodular. For conciseness, we notate both the link set and the corresponding Laplacian matrix $Q$ of the graph composed of these links by the same notation $Q$.

\begin{small}
	\begin{algorithm}[htb]
		\caption{Heuristic Greedy Method for $\mu$}\label{alg:greedy_metric}
		\begin{spacing}{1}
			\begin{algorithmic}[1]
				\Inputs {the undirected network $G$ and the integer $K$} 
				\Initialize {\strut{Set the solution $S$ as an empty link set}}
				\While {$|S| \le K$}
				\State {Compute the eigenvector $z$ of $R$ by $Q$}
				\State {$\ell_{ij} = \arg\max_{\ell\notin G} \{|z_i-z_j|\} $}
				\State {$S = S\cup \{\ell_{ij}\}$, $G = G\cup \{\ell_{ij}\}$ }		
				\EndWhile
				\State {Return $S$}
			\end{algorithmic}	
		\end{spacing}	
	\end{algorithm}
\end{small}

\begin{remark}\label{thm:submodular_lambda}
	Given the perturbation of the Laplacian matrix $Q+\Updelta Q$, the function $\mathcal{E}(\Updelta Q)$ as \eqref{equ:impact} is non-submodular with respect to the set $\Updelta Q$ of adding subgraphs (links).
\end{remark}

\textit{Proof: }
We define the Laplacian matrix for the perturbation subgraphs $I,J,U,V,W$ following the relation that $I\leq J$, $V = I+U$, $W = J+U$ and $Y = J-I$, where the symbol (+) denotes the union of link sets and (-) indicates the relative complement. 
We have the relations:
\begin{footnotesize}
	\begin{align}
	&\mathcal{E}(V)- \mathcal{E}(I)\notag\\
	&= z^T (I+U) z -{\epsilon \over \lambda_1(R)}z^T (I+U)^2 z -\left( z^T I z -{\epsilon \over \lambda_1(R)}z^T I^2 z\right) \notag\\
	&=z^T U z - {\epsilon \over \lambda_1(R)}z^T U^2 z +{\epsilon \over \lambda_1(R)}z^T (IU+UI) z \notag\\
	&=\mathcal{E}(U)-{\epsilon \over \lambda_1(R)}z^T (IU+UI) z\\
	&\mathcal{E}(W)- \mathcal{E}(J)=\mathcal{E}(U)-{\epsilon \over \lambda_1(R)}z^T (JU+UJ) z
	\end{align}
\end{footnotesize}
Thus, we can obtain that 
\begin{footnotesize}
	\begin{align}
	\left(\mathcal{E}(V)- \mathcal{E}(I)\right)-\left(\mathcal{E}(W)-\mathcal{E}(J)\right) &= {\epsilon \over \lambda_1(R)}z^T (YU+UY) z
	\end{align}
\end{footnotesize}
which is always positive only if the matrix $U$ and $Y$ are commuting.
Thus, the submodularity of the impact function $\mathcal{E}(\Updelta Q)$ is not guaranteed.  \hfill $\Box$

We then show that the impact function $\mathcal{E}(\Updelta Q)$ is approximately submodular.
\begin{lemma}\label{thm:link_number_to_lambda}
	The largest eigenvalue of the Laplacian matrix $\Updelta Q$ of a graph with $K$ links is upper bounded by $K+1$.
\end{lemma}
\textit{Proof: } The largest Laplacian eigenvalue is upper bounded by the sum of the largest nodal degree and the second largest nodal degree \cite{merris1994laplacian}.  Thus, the largest Laplacian eigenvalue of a graph with $K$ links has the upper bound $K+1$, which occurs in a star graph. \hfill $\Box$

\begin{lemma}\label{thm:approximate_modular}
	Assuming that the adding subgraph with the Laplacian matrix $\Updelta Q$ consists of $K$ links.
	The function $\mathcal{E}(\Updelta Q)$ is bounded by 
	\begin{align}
	\delta z^T\Updelta Q z \leq \mathcal{E}(\Updelta Q) \leq z^T\Updelta Q z
	\end{align}
	where the constant $\delta := 1-\frac{\epsilon}{\lambda_1(R)}(K+1)$.
\end{lemma}

\textit{Proof: } We define the constant $\delta>0$ such that
\begin{footnotesize}
\begin{align}
	\delta = \min \frac{z^T \Updelta Q z - \frac{\epsilon}{\lambda_1(R)}z^T(\Updelta Q)^2 z}{z^T \Updelta Q z}=1-\max \frac{\epsilon}{\lambda_1(R)} \frac{z^T(\Updelta Q)^2 z}{z^T \Updelta Q z}
\end{align}
\end{footnotesize}
The term $z^T(\Updelta Q)^2 z$ obey the relation
\begin{align}
	|z^T(\Updelta Q)^2 z| \leq ||\Updelta Q||\cdot|z^T \Updelta Q z|\leq \lambda_1(\Updelta Q)z^T \Updelta Q z
\end{align}
Invoking the upper-bound of $\lambda_1(\Updelta Q)\leq K+1$ in Lemma \ref{thm:link_number_to_lambda}, we can obtain that $\delta = 1-\frac{\epsilon}{\lambda_1(R)}(K+1)$ and $\delta z^T\Updelta Q z \leq \mathcal{E}(\Updelta Q)$.

The right inequality can be verified by the fact that $\frac{\epsilon}{\lambda_1(R)}z^T(\Updelta Q)^2 z\geq0$ for the positive semi-definite matrix $(\Updelta Q)^2$. \hfill $\Box$

Usually, the perturbed subgraph are a small amount of links separately located in a large network, so the eigenvalue $\lambda_1(\Updelta Q)$ is relatively small, which practically leads to the constant $0<\delta<1$. 

\begin{theorem}\label{thm:peformance_bound}
	The greedy algorithm \ref{alg:greedy_metric} for maximizing the impact function $\mathcal{E}(\Updelta Q)$ by adding $K$ links (that compose a graph with the Laplacian matrix $\Updelta Q$),  can guarantee the performance
	\begin{small}
	\begin{align} 
	\mathcal{E}(\Updelta Q) \geq \frac{1}{1+\frac{K(1-\delta^2)}{\delta^2}}\left( 1-\delta^{2K}\left(1-\frac{1}{K}\right)^K\right) \mathcal{E}_{\text{OPT}}(\Updelta Q^*)
	\end{align}
	\end{small}
if $\frac{1-\delta}{1+\delta}\leq \frac{1}{K}$, where $\mathcal{E}_{\text{OPT}}(\Updelta Q^*)$ is the optimal impact with the optimal solution set $\Updelta Q^*$.
\end{theorem}
\textit{Proof: } Lemma \ref{thm:approximate_modular} implies that the impact $\mathcal{E}(\Updelta Q)$ is a $\delta-$approximately submodular function \cite{gupta2018approximate} which approximates the (sub)modular function $g(\Updelta Q) = z^T \Updelta Q z$. 
Invoking the performance bound in Theorem 8 proposed in \cite{horel2016maximization}, we can arrive that the performance constant is $\frac{1}{1+\frac{K(1-\delta^2)}{\delta^2}}\left( 1-\delta^{2K}\left(1-\frac{1}{K}\right)^K\right)$.    \hfill $\Box$

The near optimal solution of the exact increment $\Updelta\mu^*$ can be approached by Algorithm  \ref{alg:greedy_metric} in some networks.
Theorem \ref{thm:peformance_bound} shows that the performance of the greedy algorithm is related to the topological property of the original network indicated by the constant $\theta$, i.e., a smaller $\theta=\frac{\epsilon}{\lambda_1(R)}$ leading a larger $\delta$ yields a larger performance constant.  
Specifically, the impact $\mathcal{E}(\Updelta Q)$ approximates to the modular(additive) function $g(\Updelta Q) $ for a smaller $\theta$, which implies a better performance for the greedy algorithm. 
The performance constant also degrades with the increasing number of adding links~$K$. 
In addition, the computational complexity of Algorithm \ref{alg:org_greedy} is $\mathcal{O}(K\bar{L}N^2)$, while the computational complexity of Algorithm \ref{alg:greedy_metric} reduces to $\mathcal{O}(K(\bar{L}+N^2))$, where $\bar{L} = {N\choose 2}-L$ is the number of links in the complement of $G$.


\section{Generalized algebraic connectivity in directed networks}
\begin{figure}[tb]\centering
	\includegraphics[width=9cm]{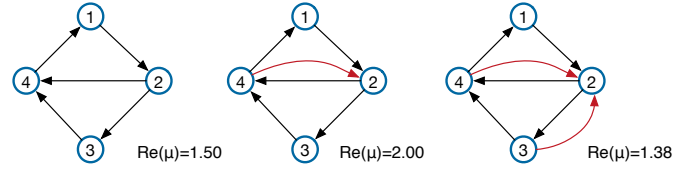}
	\caption{Example of the non-monotonicity of the generalized algebraic connectivity $\Re(\mu)$ by adding links.}
	\label{fig:nonmono}
\end{figure}

This section extends our investigation to directed networks.
We define $\Re(\mu)$ as the second smallest real part among all the eigenvalues of the generalized Laplacian matrix $Q =\Delta_{in}-A^T$ in directed networks without self-loops. 
The properties of the generalized algebraic connectivity $\Re(\mu)$ in directed networks is more complicated than those in undirected networks. For example, the generalized algebraic connectivity $\Re(\mu)$ is even not monotonic with the number of adding links, and an additional directed link could decrease the generalized algebraic connectivity (see an example in Figure \ref{fig:nonmono}). 

The algebraic connectivity in a directed network for a real perturbation of Laplacian matrix $\Updelta Q$ is given by $\Updelta\mu = \frac{y^H\Updelta Q x}{y^Hx}+ \mathcal{O}(\Vert\Updelta Q\Vert_F^2 )$, where $x$ and $y$ are right- and left- eigenvectors belonging to $\mu$ normalized so that $\Vert x\Vert_2=\Vert y\Vert_2=1$ and $y^Hx =|y^Hx|$. For the complex eigenvalue and eigenvectors, we have the perturbation formula
\begin{align}\label{equ:real_perturbation_expansion}
\Re(\Updelta\mu) = \frac{y_R^T\Updelta Q x_R+y_I^T\Updelta Q x_I}{y_R^T x_R+y_I^T x_I}+ \mathcal{O}(\Vert\Updelta Q\Vert_F^2 )	
\end{align}
where $x=x_R+ix_I$ and $y=y_R+iy_I$. Although the perturbation formula \eqref{equ:real_perturbation_expansion} provides an estimation of $\Re(\Updelta\mu)$, the computational cost of $x$ and $y$ is very high in large networks. Moreover, the change of the real part of a determined eigenvalue $\Re(\Updelta\mu)$ could not equal the change of the generalized algebraic connectivity $\Updelta\Re(\mu)$. 

\subsection{Bounds of the generalized algebraic connectivity perturbation $\Updelta\Re(\mu)$}
Unfortunately, the complex algebraic connectivity $\mu$ of a non-Hermitian Laplacian matrix does not obey the power iteration method, which leads that the generalized algebraic connectivity $\Re(\mu)$ perturbation difficult to derive in a similar form with \eqref{equ:mu_perturbation}. 
Thus, we consider a lower-bound and an upper-bound of the generalized algebraic connectivity $\Updelta\Re(\mu)$ under topological perturbation.

\begin{lemma}\label{thm:positive_define}	
	Given the generalized Laplacian matrix $Q =\Delta_{in}-A^T$, the matrix $Q+Q^T$ is positive semi-definite if the network is strongly connected.
\end{lemma}

\textit{Proof: } By the Gershgorin Disk Theorem \cite{van2010graph}, all the eigenvalues are located in the union of the disks centered at $Q_{jj}+ Q^T_{jj}= \sum_{i=1}^N a_{ij} + \sum_{i=1}^N a_{ji}$ with the radius $r_j = \sum_{i=1}^N a_{ij}+\sum_{i=1}^N a_{ji}$ for $j\in\mathcal{N}$. Thus, all the eigenvalues are located in the right plane, which yields that the symmetric matrix $Q+Q^T$ is positive semi-definite. \hfill $\Box$

\begin{theorem}\label{thm:generalize_lower_bound}
	Given a network with the generalized Laplacian matrix $Q$, the increment of the generalized algebraic connectivity $\Updelta\Re(\mu(\Updelta Q))$ by adding a subgraph with the generalized Laplacian matrix $\Updelta Q$ can be lower-bounded by an approximately submodular function, i.e.,
	\begin{align}\label{equ:lower_bound_thm}
	\Updelta\Re(\mu(\Updelta Q)) \geq \frac{1}{2}z^T\Updelta Q^* z  -\frac{\epsilon}{4\lambda_1(H)}z^T(\Updelta Q^*)^2 z
	\end{align}
	where $\Updelta Q^* = \Updelta Q+\Updelta Q^T$, $z$ is the principle eigenvector of $H = \frac{1}{2}(R+R^T)$, $R$ and $\epsilon$ are defined in \eqref{equ:defineR}.
\end{theorem}
\textit{Proof: } According to Section 3.1, we have the shift relation that $\lambda_i(Q) = \frac{1}{\epsilon}(1-\lambda_{N-i}(R))$ for the matrix $R = I_N - \epsilon Q-\frac{1}{N}J$. Similarly, we can obtain that 
\begin{align}
\Re(\mu(Q)) = \frac{1}{\epsilon}\left(1-\max_{i\geq 2} \Re(\lambda_{i-1}(R))\right)
\end{align}
We then introduce the Bendixson theorem \cite{kato2013perturbation} that: let $H=\frac{1}{2}(R+R^T)$ be the Hermitian part of $R$, we have $\lambda_N(H) \leq \Re(\lambda(R)) \leq \lambda_1(H)$. 
For the perturbation $Q+\Updelta Q$, we have $\Updelta R = -\epsilon \Updelta Q$, which yields $\Updelta H = -\frac{\epsilon}{2} (\Updelta Q+\Updelta Q^T) = -\frac{\epsilon}{2}\Updelta Q^*$, where $\Updelta Q^*$ is positive semi-definite due to Lemma \ref{thm:positive_define}.
Invoking the eigenvalue perturbation in \eqref{equ:mu_perturbation}, we arrive at
\begin{small}
\begin{align}\label{equ:real_upperbound}
\Re(\lambda_{i-1}(R+\Updelta R)) &\leq \lambda_{1}(H+\Updelta H) \notag\\
&\approx \lambda_{1}(H)-\frac{\epsilon}{2}z^T\Updelta Q^* z  +\frac{\epsilon^2}{4\lambda_1(H)}z^T(\Updelta Q^*)^2 z	
\end{align}
\end{small}
where $z$ is the left eigenvector of the matrix $H$ corresponding to the largest eigenvalue $\lambda_1(H)$ normalized by $\Vert z\Vert_2 =1$. 
We can obtain the lower bound of the increment of the general algebraic connectivity $\Updelta\Re(\mu(\Updelta Q))$ follows
\begin{small}
\begin{align}
\Updelta\Re(\mu(\Updelta Q)) &= \Re(\mu(Q+\Updelta Q))-\Re(\mu(Q)) \notag\\
&\geq \frac{1}{2}z^T\Updelta Q^* z  -\frac{\epsilon}{4\lambda_1(H)}z^T(\Updelta Q^*)^2 z 
\end{align}
\end{small}
Since the matrix $\Updelta Q^*$ is positive definite and Hermitian, the function on the right side of \eqref{equ:lower_bound_thm} is an approximately submodular function (similar with Lemma \ref{thm:approximate_modular}), which completes the proof. \hfill $\Box$.

We also present an upper-bound of the generalized algebraic connectivity perturbation $\Updelta\Re(\mu)$ based on Lemma~\ref{thm:exponential_convert}.
\begin{lemma}\cite{asadi2017expected}\label{thm:exponential_convert}
	Defining the matrix for a directed network as $R = e^{I_N-\alpha Q}-e x_N x_N^T$, where $x_N$ is the right eigenvector of the Laplacian matrix $Q$ associated with its zero eigenvalue, and $0<\alpha<d_{in}^{\max}$, the generalized algebraic connectivity holds $\Re(\mu(Q)) = \frac{1}{\alpha}(1-\log(\max_{i\in N} |\lambda_i(R)|))$.
\end{lemma}
According to Lemma \ref{thm:exponential_convert}, computing the generalized algebraic connectivity $\Re(\mu)$ can be reduced to computing the maximum absolute eigenvalue of the matrix $R$ by an exponential operator. Then, the generalized algebraic connectivity $\Re(\mu)$ in large scale networks can be computed by a general power iteration method by Krylov subspace \cite{ruhe1994rational}. 

\begin{theorem}\label{thm:exponential_upperbound}
	The increment of the generalized algebraic connectivity $\Updelta\Re(\mu(\Updelta Q))$ by adding a subgraph with the generalized Laplacian matrix $\Updelta Q$ can be upper-bounded by the function
	\begin{align}\label{equ:upper_bound_thm}
	\Updelta\Re(\mu(\Updelta Q)) \leq \frac{1}{\alpha}\log\frac{\widehat{\lambda}(R)}{\widehat{\lambda}(R)-\kappa(Z)\Vert \alpha\Updelta Qe^{I_N-\alpha Q}\Vert_p}
	\end{align}
	where $\alpha$ and $R$ are defined in Lemma \ref{thm:exponential_convert}, $\kappa(Z)$ is the condition number of the matrix $Z$ composed by the eigenvector of $R$, and $\widehat{\lambda}(R) = \max_{i\in N} |\lambda_i(R)|$.
\end{theorem}
\textit{Proof: }We have the perturbation of the matrix
\begin{footnotesize}
\begin{align*}
R+\Updelta R = e^{I_N-\alpha(Q+\Updelta Q)}-e x_N x_N^T\approx(I_N-\alpha\Updelta Q)e^{I_N-\alpha Q}-e x_N x_N^T
\end{align*}
\end{footnotesize}
which yields $\Updelta R \approx -\alpha\Updelta Qe^{I_N-\alpha Q}$ for a small $\alpha$. According to Lemma \ref{thm:exponential_convert}, the problem of maximizing the general algebraic connectivity is equal to minimizing the largest absolute eigenvalue of $R$. Invoking the Bauer-Fike theorem \cite{bauer1960norms} and the reverse triangle inequality, we can obtain that
\begin{footnotesize}
\begin{align}
|\lambda(R)|-|\lambda(R+\Updelta R)| 
\leq |\lambda(R)-\lambda(R+\Updelta R)|\leq \kappa(Z)\Vert \Updelta R \Vert_p
\end{align}
\end{footnotesize}
where $\kappa(Z)= \Vert Z\Vert_p \Vert Z^{-1}\Vert_p$ is the condition number \cite{belsley2005regression} of the matrix $Z$ composed by the eigenvector of $R$.

Denoting the maximum absolute eigenvalue $\max_{i\in N} |\lambda_i(R)|$ by $\widehat{\lambda}(R)$, we have
\begin{small}
\begin{align}
\Updelta\Re(\mu(\Updelta Q)) &= \frac{1}{\alpha}(1-\log\widehat{\lambda}(R+\Updelta R)) - \frac{1}{\alpha}(1-\log\widehat{\lambda}(R)) \notag\\
&\leq \frac{1}{\alpha}\log\frac{\widehat{\lambda}(R)}{\widehat{\lambda}(R)-\kappa(X)\Vert \Updelta R\Vert_p} \notag\\
&=\frac{1}{\alpha}\log\frac{\widehat{\lambda}(R)}{\widehat{\lambda}(R)-\kappa(X)\Vert \alpha\Updelta Qe^{I_N-\alpha Q}\Vert_p}
\end{align}
\end{small}
which completes the proof. \hfill $\Box$

\subsection{Maximize the generalized algebraic connectivity by adding links}
We then consider the problem of maximizing the generalized algebraic connectivity by adding $K$ links in directed networks. Assuming that we apply the greedy strategy to adding the links one by one, a simple and intuitive method to select the link in each iteration follows two steps: 1. obtain the best undirected link regarding the network as undirected; 2. randomly determine the direction of this link. However, the performance of this method is not guaranteed since the effect of the link direction is ignored.
The bounds of the perturbation indicate the impact of topological perturbation in two different aspects: the lower bound \eqref{equ:lower_bound_thm} implies that the increment of the generalized algebraic connectivity is related to the connections between nodes with different eigenvector centrality, while the exponential norm of the perturbation $\Vert\Updelta R\Vert$ influences the upper bound \eqref{equ:upper_bound_thm}.

Further, Theorem \ref{thm:generalize_lower_bound} demonstrates that the maximal impact of an individual additional link on the lower bound of the increment of the general algebraic connectivity $\Updelta\Re(\mu(\Updelta Q))$ can be measured by
\begin{align}\label{equ:lower_bound}
\underline{\Omega}(\ell_{ij}) = \frac{1}{2}z_j(z_j-z_i)-\frac{\epsilon}{4\lambda_1(H)}\left((2z_j-z_i)^2+z_j^2\right)
\end{align}
which is different from the impact $\Omega(\ell_{ij}) = |z_i-z_j|$ in undirected networks.
The dominating first term $\frac{1}{2}z_j(z_j-z_i)$ in \eqref{equ:lower_bound} implies that the best link $\ell_{ij}$ tends to be located in two different communities to maximize the difference $z_j-z_i$.
Meanwhile, for the directed links between the same two nodes, i.e., $\ell_{ij}$ and $\ell_{ji}$, the target of the additional link tends to be the node with a higher value $z_j$ to increase the generalized algebraic connectivity. 
Theorem \ref{thm:exponential_upperbound} demonstrates that the impact of the subgraph depends on $\Vert\Updelta R\Vert_p$. Defining the matrix $W = e^{I_N-\epsilon Q}\approx 2I_N-\epsilon Q+\frac{1}{2}(I_N-\epsilon Q)^2$ and applying norm-1 for computational simplicity, we can obtain another metric to measure the importance of an individual link $\ell_{ij}$ on the generalized algebraic connectivity, i.e.,
\begin{align} \label{equ:upper_bound}
\overline{\Omega}(\ell_{ij}) = \Vert\Updelta Qe^{I_N-\epsilon Q}\Vert_1 = \max_{k\in N}|W_{ik}-W_{jk}|
\end{align}	

We assume that the impact of an adding link on the bounds of $\Updelta\Re(\mu)$ has similar behaviors with the exact generalized algebraic connectivity.
By adding the link with the largest impact in each iteration, we heuristically propose the greedy method as Algorithm \ref{generalized algebraic connectivity} to approach the near-optimal solution.

\begin{small}
	\begin{algorithm}[htb]
		\caption{Heuristic Greedy Method for $\Re(\mu)$ }\label{generalized algebraic connectivity}
		\begin{spacing}{1}
			\begin{algorithmic}[1]
				\Inputs {the undirected network $G$ and the integer $k$} 
				\Initialize {\strut{Set the solution $S$ as an empty link set}}
				\While {$|S| \le k$}
				\State {Compute matrix $R$ by the Laplacian matrix $Q$}
				\State {$\ell_{ij} = \arg\max_{\ell\notin G} \underline{\Omega}(\ell_{ij}) $ or $\ell_{ij} = \arg\max_{\ell\notin G} \overline{\Omega}(\ell_{ij})$}
				\State {$S = S\cup \{\ell_{ij}\}$, $G = G\cup \{\ell_{ij}\}$ }		
				\EndWhile
				\State {Return $S$}
			\end{algorithmic}	
		\end{spacing}	
	\end{algorithm}
\end{small}

\section{Numerical evaluations}
In this section, we present some numerical tests to evaluate the performance of the perturbation formula and our proposed methods.
The topological properties of the investigated networks are summarized in Table \ref{tab:netinfo}. We extract the giant component from the undirected networks and extract the largest strongly connected component from the directed networks.
\begin{footnotesize}
\begin{table}[!htb]
	\begin{footnotesize}
	\centering
	\begin{tabular}{l|lllll}
		\hline
		& $N$ & $L$ & $\lambda_1$ & $\Re(\mu)$ & Type \\ \hline
		Karate \cite{zachary1977information}         & 34  & 78 & 6.73 & 0.469 & undirected \\
		Les Mis\'{e}rables \cite{knuth1993stanford}         & 77  & 254 & 12.00 & 0.205 & undirected \\
		NetScience \cite{newman2006finding}   & 379 & 914 & 10.38 & 0.015 & undirected \\ 
		Illinois friendship \cite{coleman1964introduction}     & 67 & 359 & 5.780 & 0.115 & directed \\ 
		Berlin traffic \cite{jahn2005system}    & 216 & 514 & 3.35 & 0.022 & directed \\ 
		Neural network \cite{white1986structure}    & 239 & 1912 & 9.15 & 0.364 & directed \\ 
		\hline
	\end{tabular}
	\caption{The topological properties of the giant component or the largest strongly connected component of several experimental networks}
	\label{tab:netinfo}
\end{footnotesize}
\end{table}
\end{footnotesize}

\subsection{Undirected networks}
\begin{figure}[tp]
	\centering
	\subfloat[Circle with 20 nodes]
	{\includegraphics[width=4.4cm]{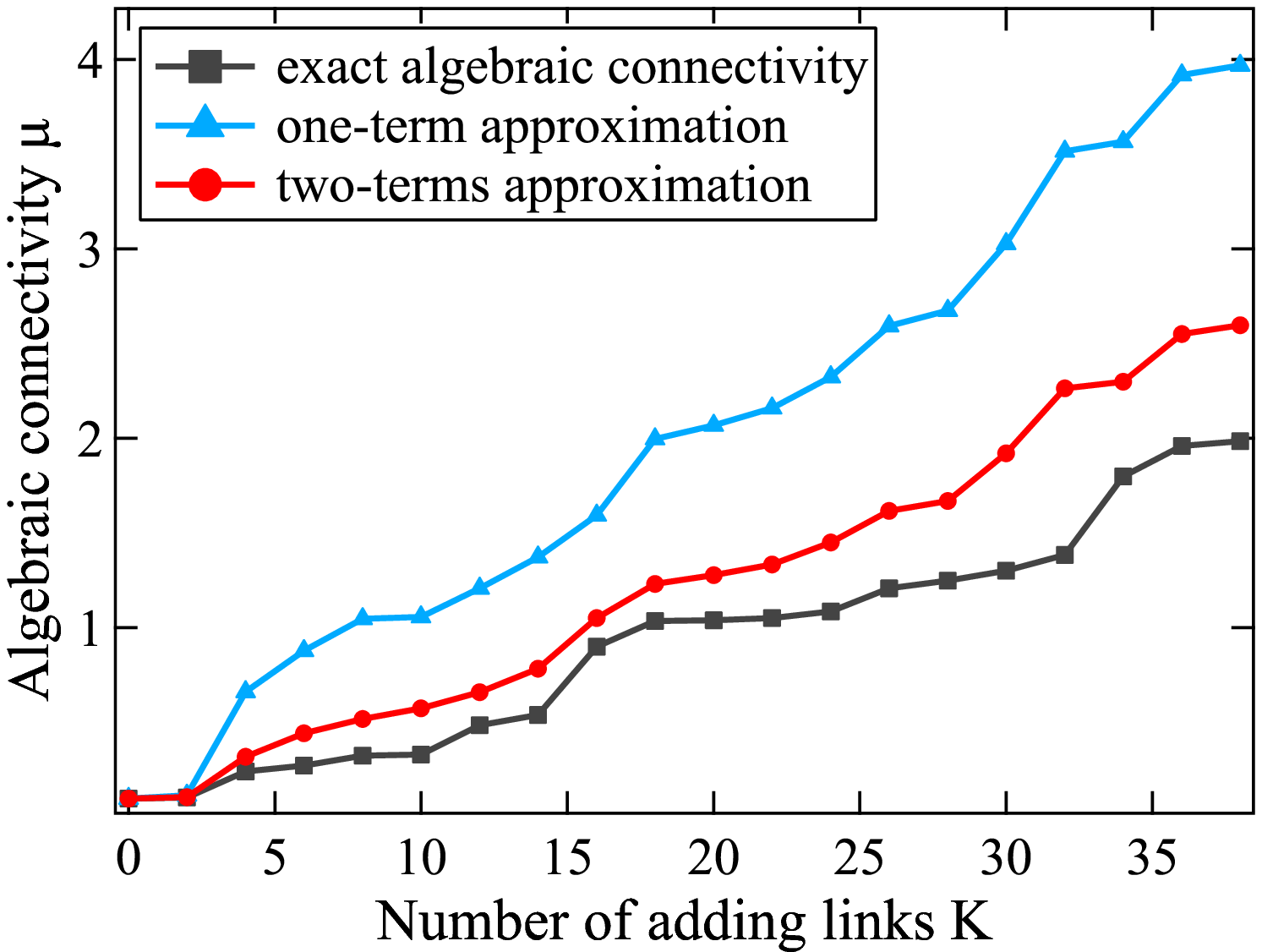}} \hfil
	\subfloat[Karate]
	{\includegraphics[width=4.4cm]{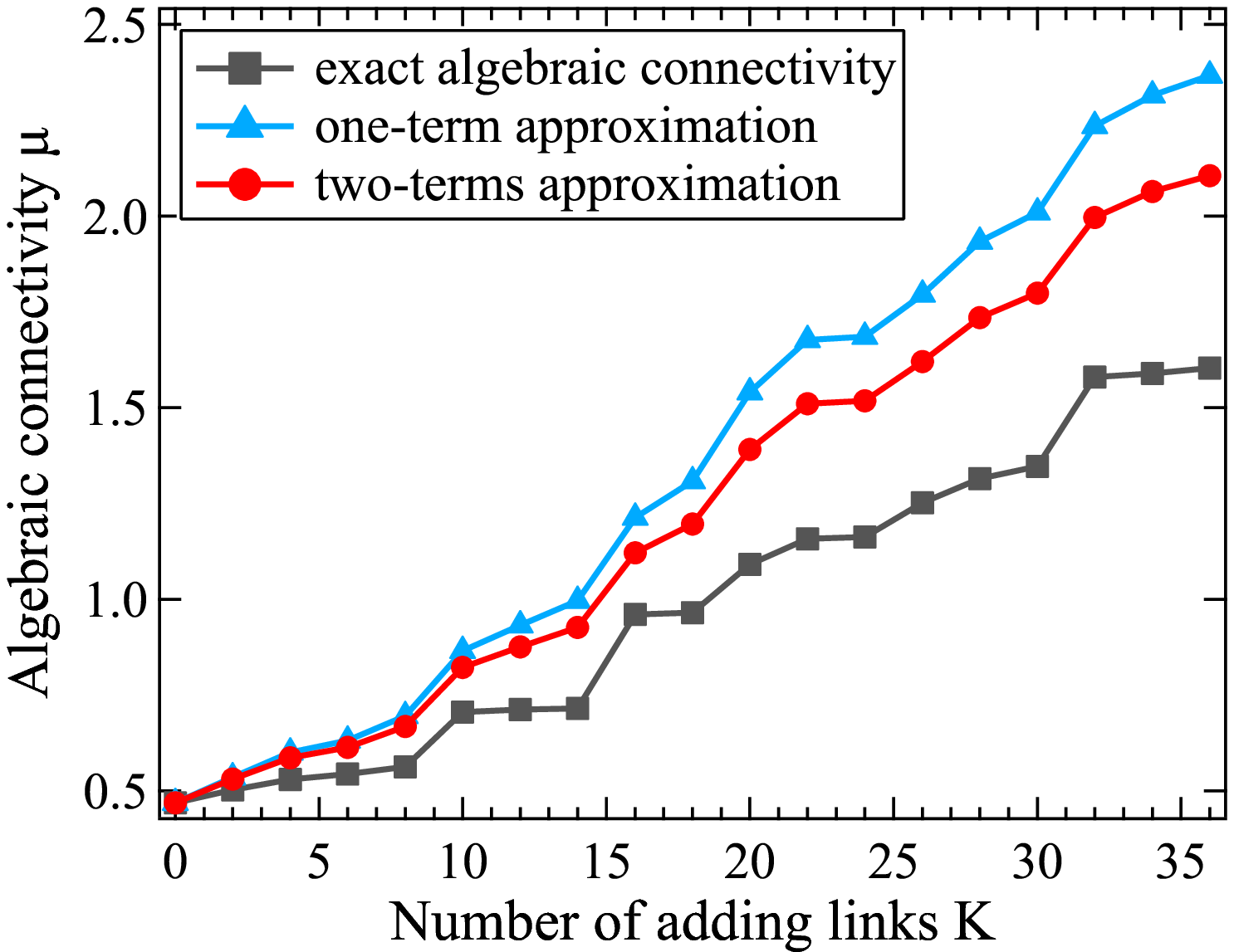}}
	\caption{The functions of (i) the exact algebraic connectivity, (ii) the approximation with one term $\widetilde{\mu} \approx \mu + z^T \Updelta Q z$ and (iii) the approximation \eqref{equ:mu_perturbation} with two terms with respected to the number of adding links $K$. }
	\label{fig:perturbation_equation}
\end{figure}

We first evaluate the perturbation formula \eqref{equ:mu_perturbation} of the algebraic connectivity by continuously adding links in undirected networks. In each step, we randomly select one node and randomly adding two new links incident to this node. 
Figure \ref{fig:perturbation_equation} shows the exact algebraic connectivity as a function of the number of the adding links $K$, which is compared with the perturbation formula with one term, i.e., $\widetilde{\mu} \approx \mu + z^T \Updelta Q z$ and the perturbation formula \eqref{equ:mu_perturbation} with two terms $\widetilde{\mu} \approx \mu + \Updelta\mu$.
We observe that the algebraic connectivity estimated by the perturbation formula \eqref{equ:mu_perturbation} usually provides an upper-bound of the exact algebraic connectivity $\mu$, but deviates more from the exact algebraic connectivity $\mu$ with the increasing number of adding links, which implies that the performance of the perturbation formula could degrade with the increasing size of the adding subgraphs. 
The proposed perturbation formula \eqref{equ:mu_perturbation} can provide a better estimation for the exact algebraic connectivity $\mu$, especially for the spare network with a small maximum degree (e.g., circle networks). 
The proposed approximation \eqref{equ:mu_perturbation} exhibits a similar behavior with the perturbation with one term in a dense network since the small $\epsilon < \frac{1}{d_{\max}}$ of \eqref{equ:mu_perturbation} reduces the impact of the second term.

We then compare the performance of different strategies to maximize the algebraic connectivity, which includes the strategies for selecting a link in each iteration: 
\begin{itemize}
	\item[(1)] Alg. 1: selecting the link $\ell_{ij}$ to maximize $\mu$ greedy; 
	\item[(2)] Alg. 2: selecting the link $\ell_{ij}$ with $\max\{|z_i-z_j|\}$; 
	\item[(3)] Selecting the link $\ell_{ij}$ between the nodes with the smallest degree product $d_id_j$; 
	\item[(4)] Selecting the link $\ell_{ij}$ between the nodes with the smallest eigenvector centrality product; 
	\item[(5)] Selecting the link $\ell_{ij}$ between the nodes with the smallest node betweenness product.
\end{itemize}

\begin{figure*}[htp]
	\centering 
	\subfloat[Karate]
	{\includegraphics[width=5.8cm]{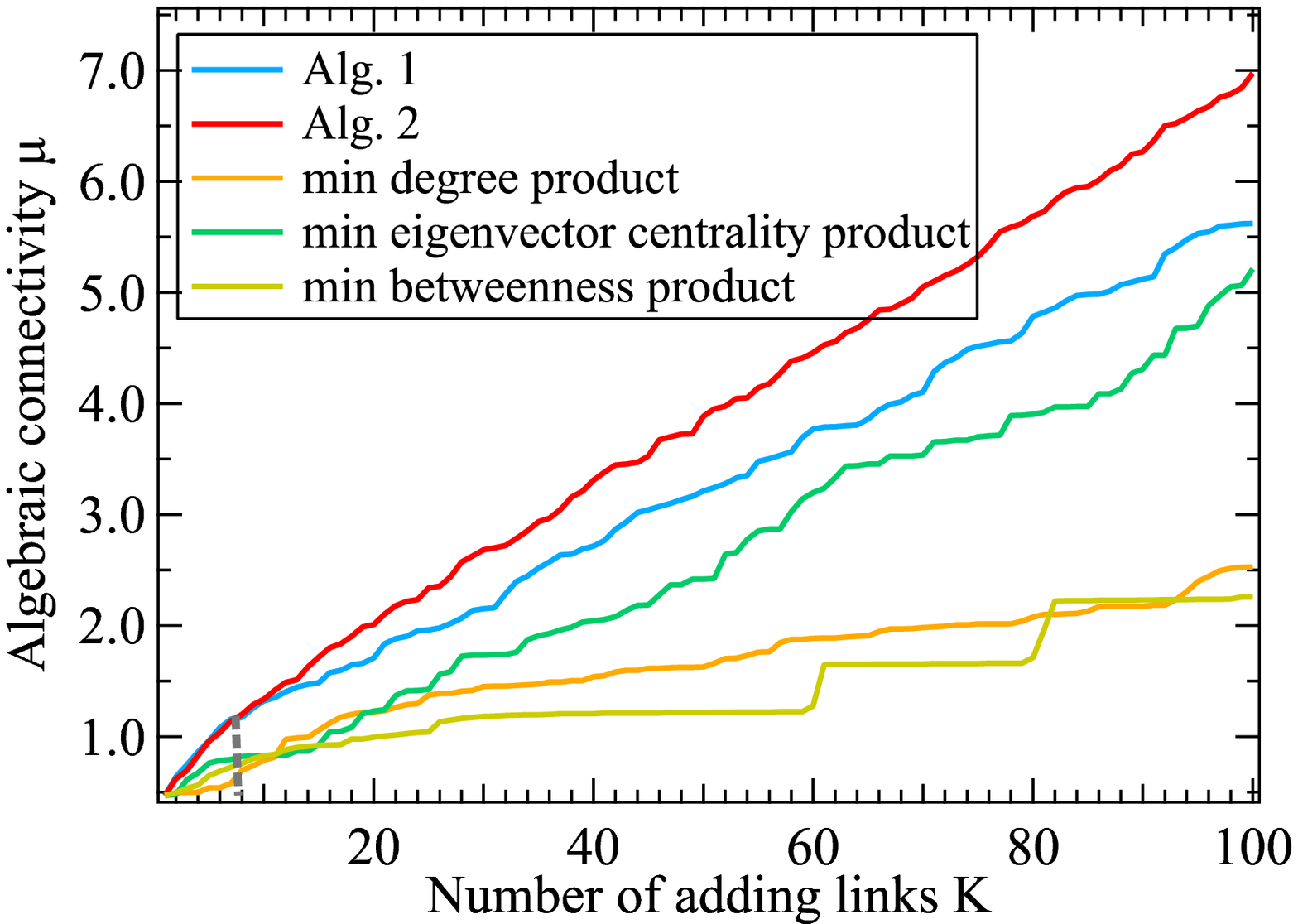}}
	\subfloat[Les Mis\'{e}rables]
	{\includegraphics[width=5.8cm]{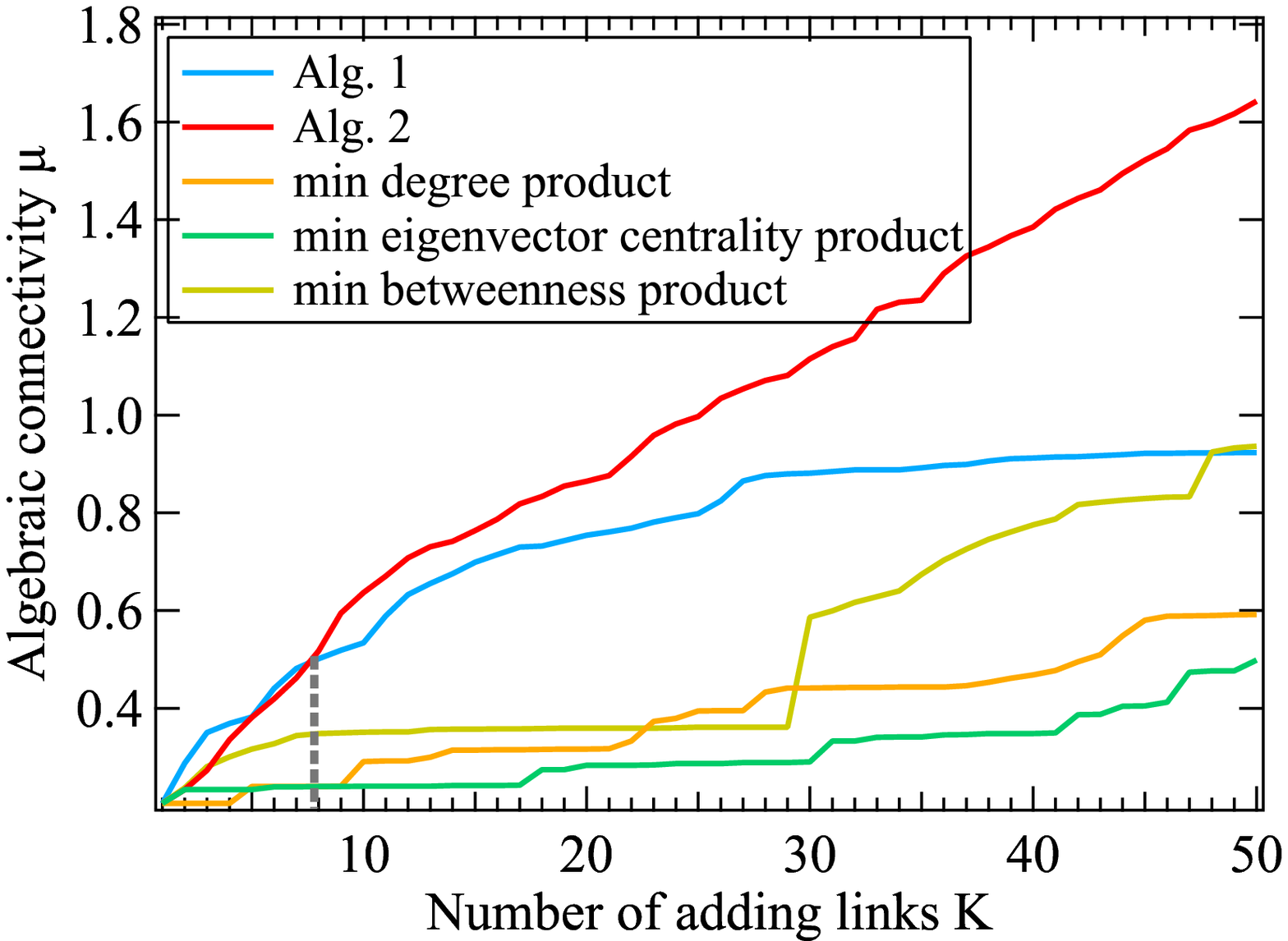}}	
	\subfloat[Netscience]
	{\includegraphics[width=5.8cm]{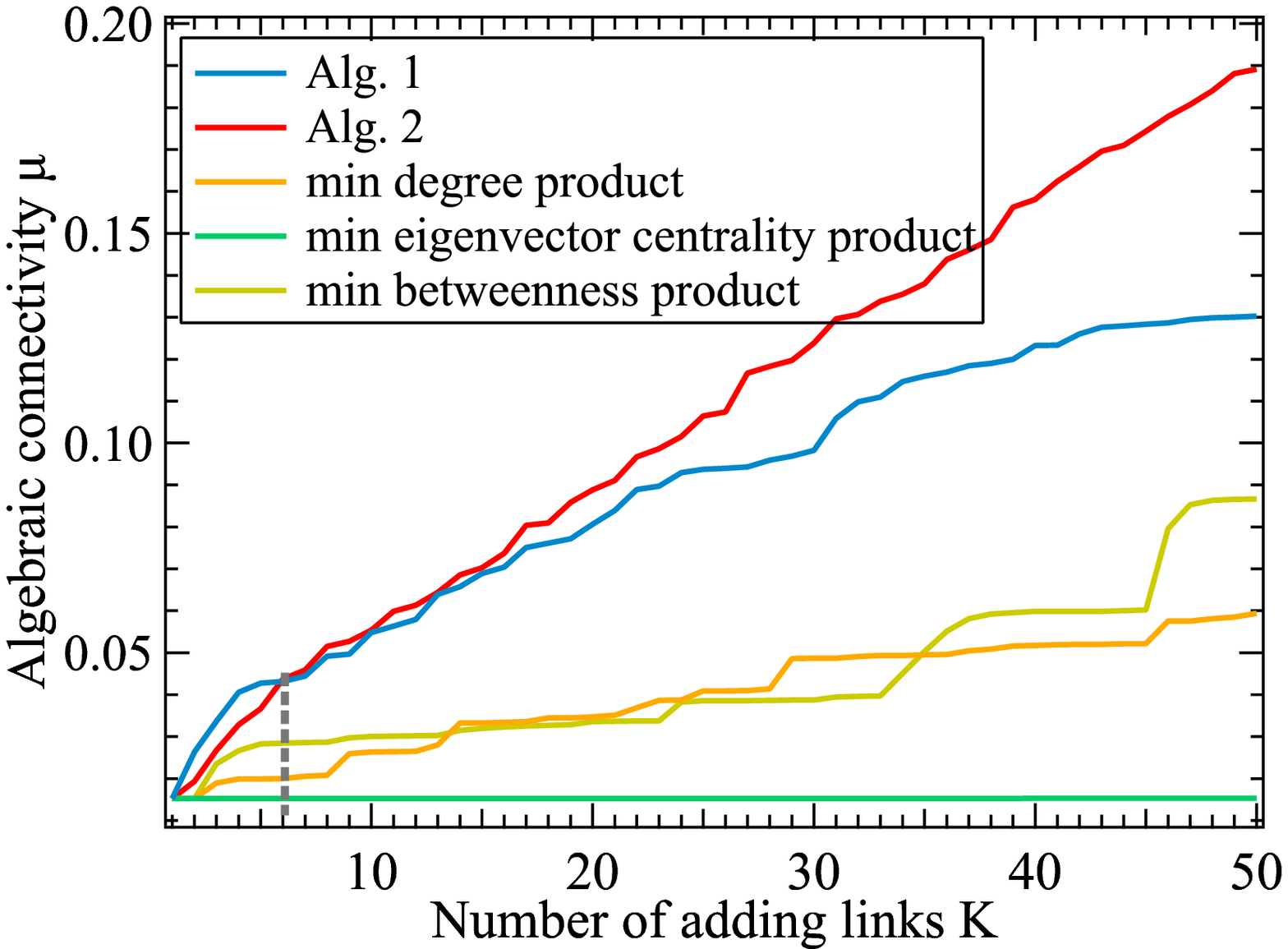}}
	\caption{The exact algebraic connectivity $\mu$ as a function of the number of adding links $K$ via different greedy strategies in three empirical undirected networks including Karate, Les Mis\'{e}rables and Netscience.}
	\label{fig:undir_performance}
\end{figure*}

Figure \ref{fig:undir_performance} shows the exact algebraic connectivity $\mu$ as a function of the number of adding links $K$ via different greedy strategies in three empirical undirected networks.
Figure \ref{fig:undir_performance} shows that the algebraic connectivity $\mu$ via Algorithm \ref{alg:org_greedy} exhibits a concave-like function with respect to the number of adding links $K$, and has the best performance for a small number of adding links (e.g., $K\approx5$). 
The proposed heuristic method Algorithm \ref{alg:greedy_metric} approaches the performance of Algorithm \ref{alg:org_greedy} best, and even outperforms Algorithm \ref{alg:org_greedy} when the number of adding links $K$ is large.
We also observe that the other heuristic greedy strategies based on the topological metrics (e.g., degree, eigenvector centrality and betweenness) usually cannot guarantee the performance for maximizing the algebraic connectivity.

\subsection{Directed networks}
\begin{figure*}[!htp]
	\centering
	\subfloat[Karate]{\includegraphics[width=6.5cm]{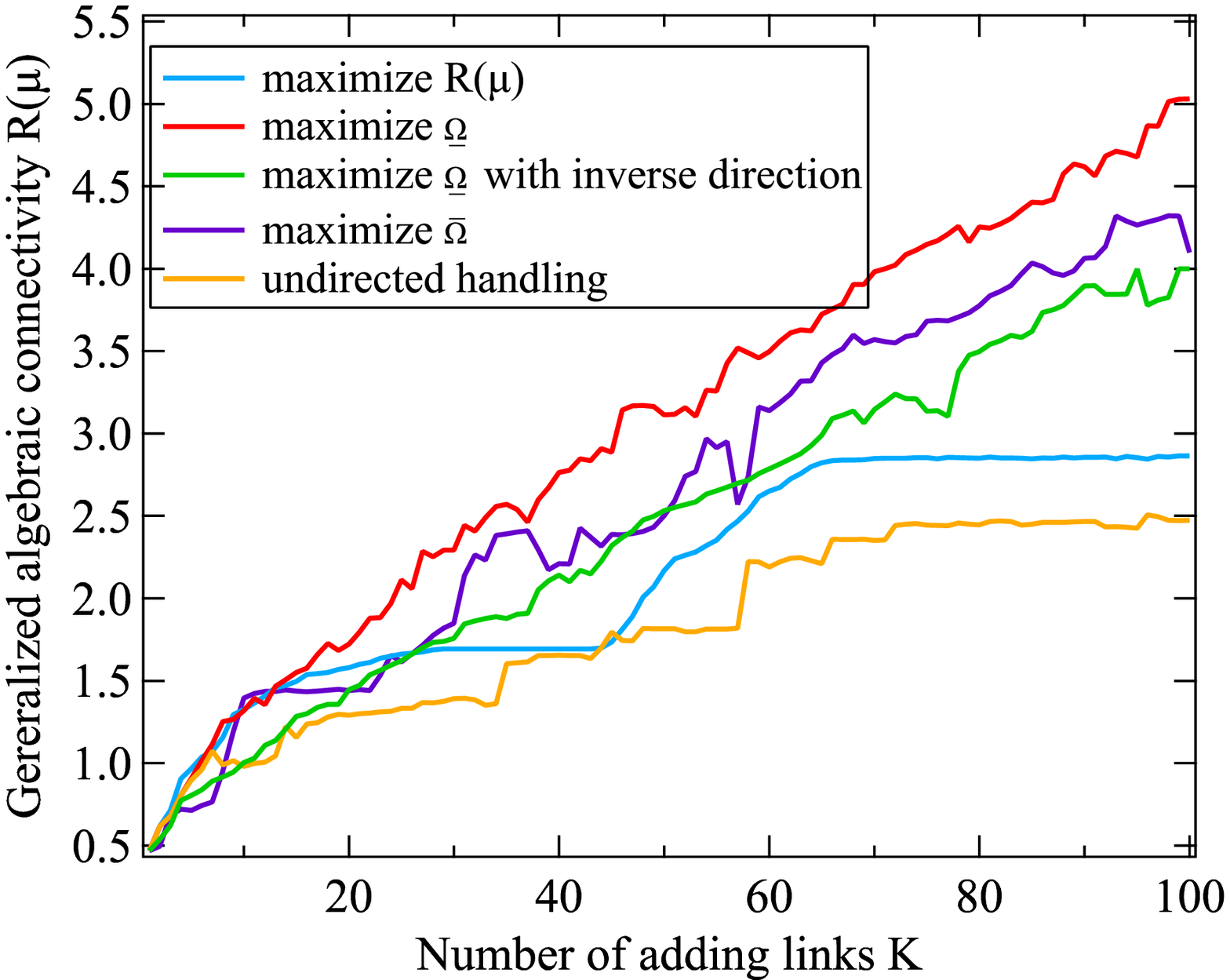}}\hfil
	\subfloat[Neural network]{\includegraphics[width=6.5cm]{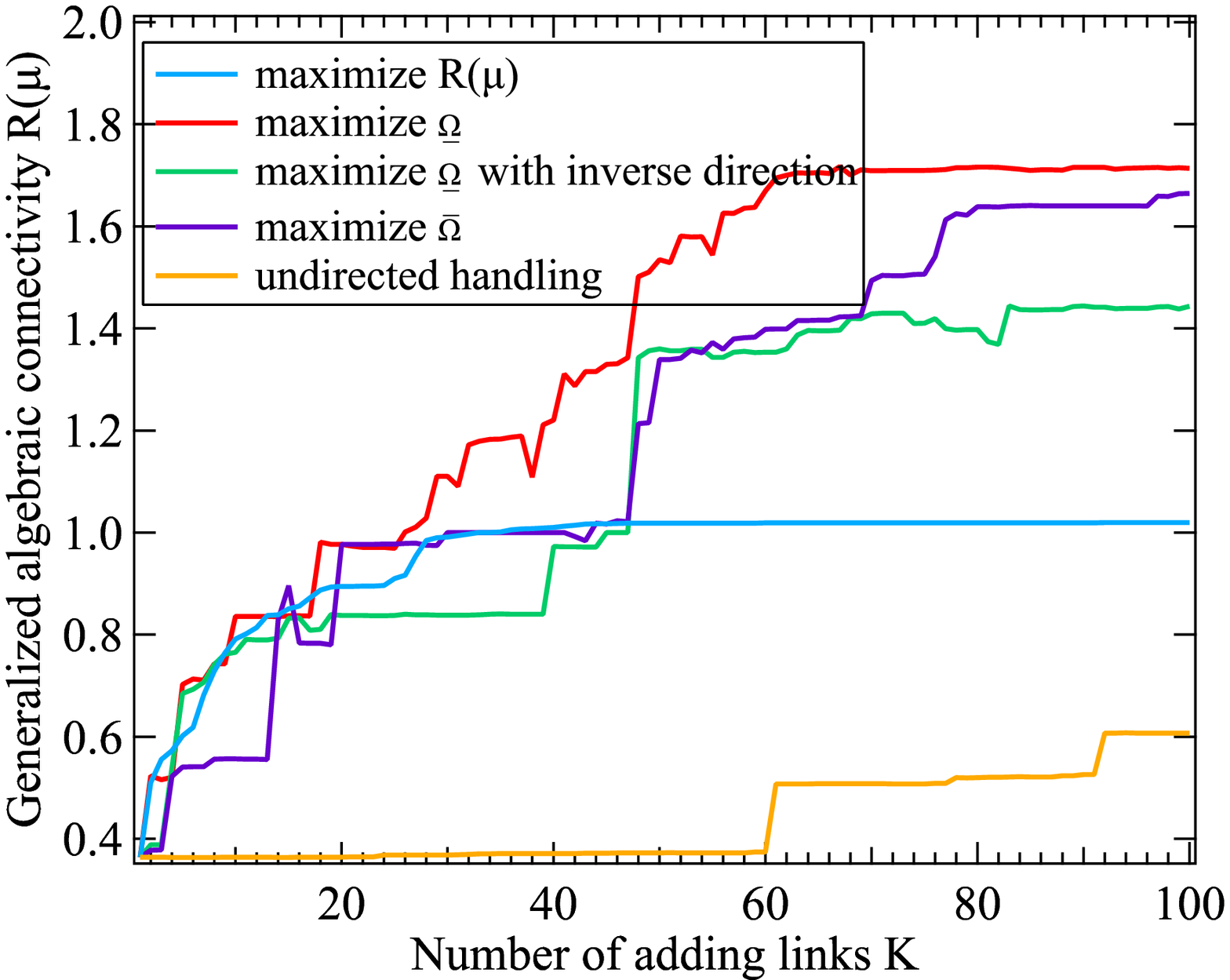}} \\	
	\subfloat[Illinois friendship]
	{\includegraphics[width=6.5cm]{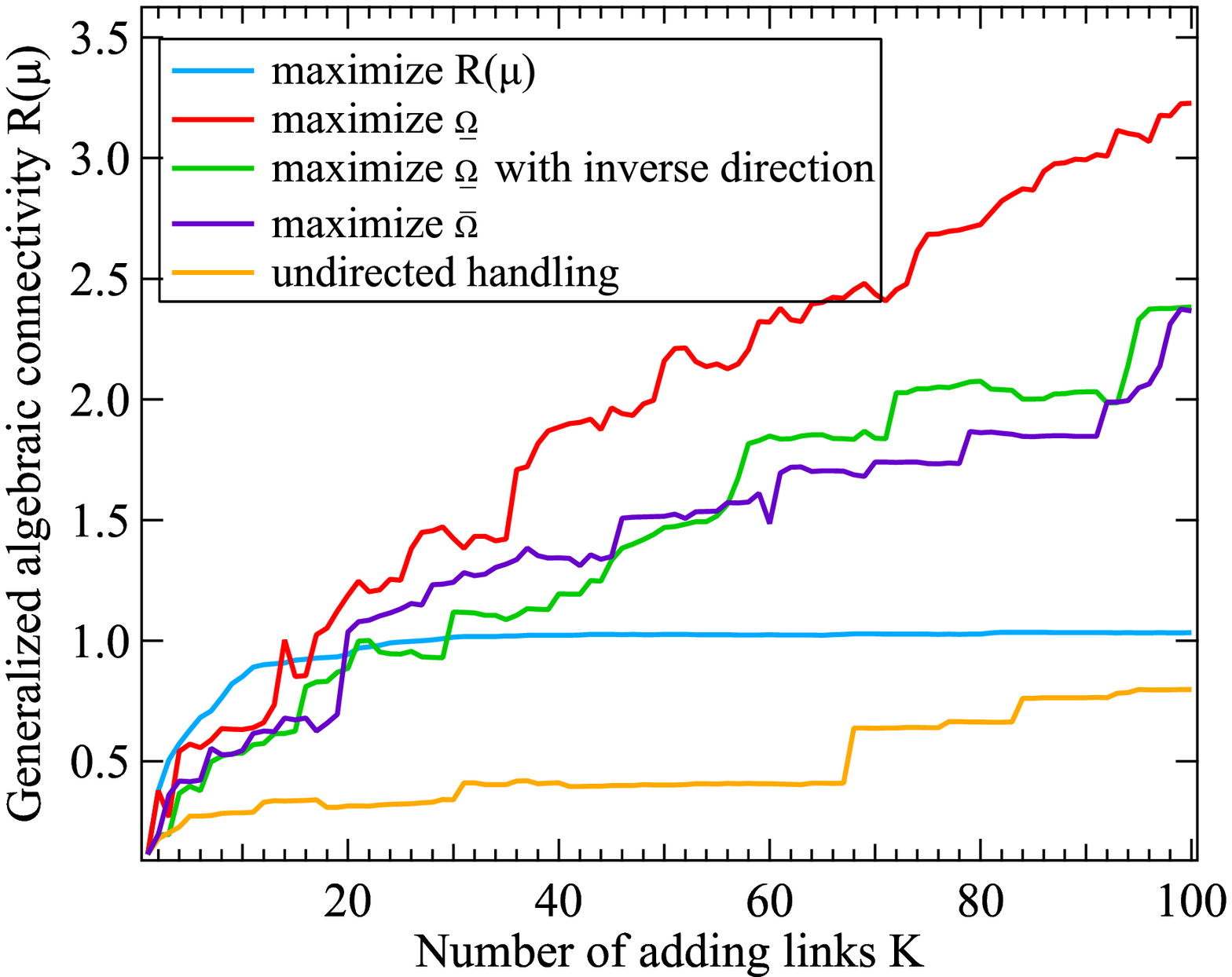}}\hfil
	\subfloat[Berlin traffic network]{\includegraphics[width=6.5cm]{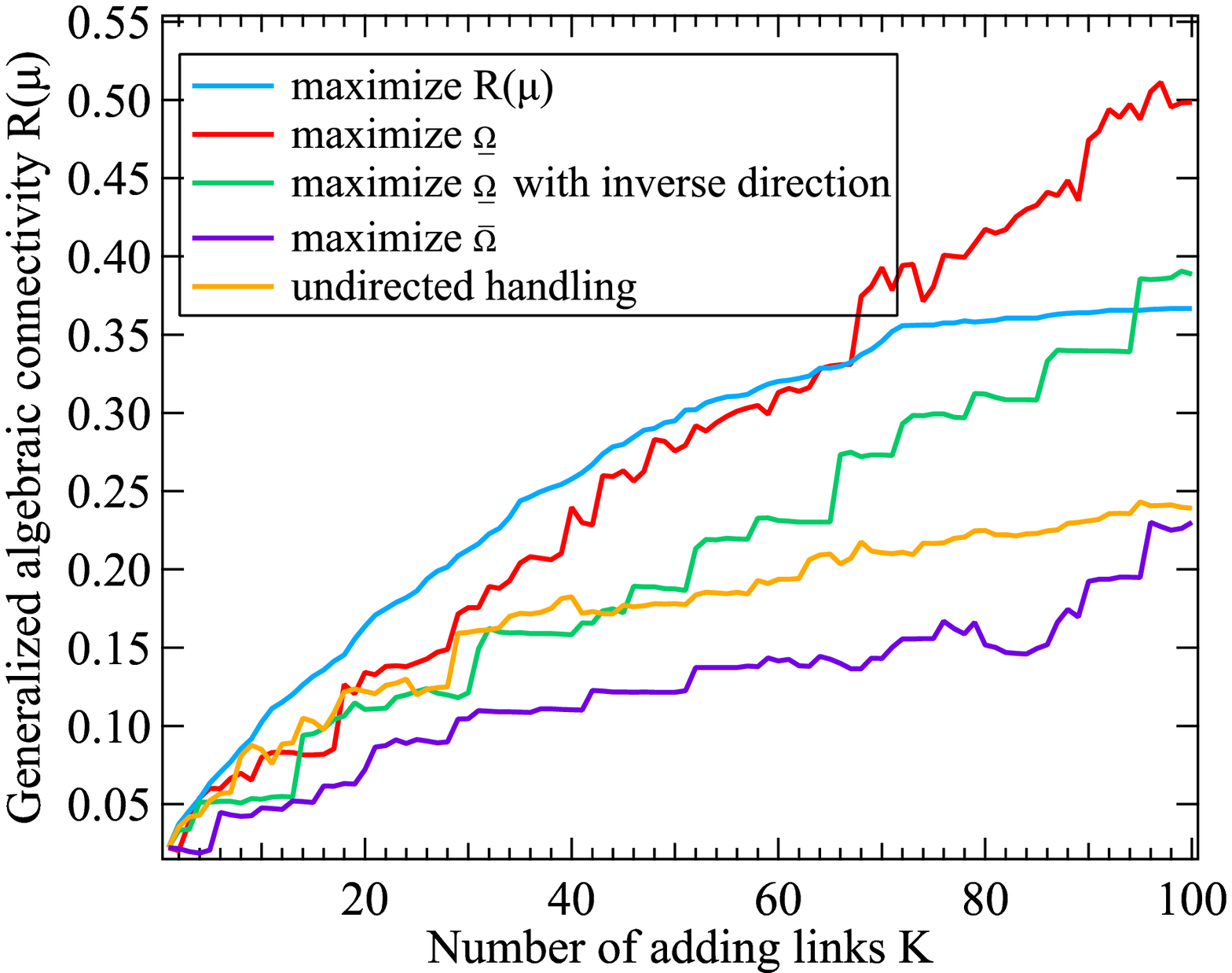}}	
	\caption{The exact generalized algebraic connectivity $\Re(\mu)$ as a function of the number of adding links $K$ via different greedy strategies in four empirical undirected networks including bi-directed Karate, Neural network, Illinois friendship and Berlin traffic network.}
	\label{fig:performance_dir}
\end{figure*}

We further evaluate the performance of the heuristic metrics proposed in Section 4.3 in directed networks. We consider the following metrics to select the individual link in each iteration in the greedy strategy:
\begin{itemize}
	\item[(1)] Maximize $\Re(\mu)$: the link that maximizes the generalized algebraic connectivity $\Re(\mu)$;
	\item[(2)] Maximize $\underline{\Omega}(\ell_{ij})$: the link with the maximum $\underline{\Omega}(\ell_{ij})$;
	\item[(3)] Maximize $\underline{\Omega}(\ell_{ij})$ with inverse direction: the inverse directed link of the link via (2); 
	\item[(4)] Maximize $\overline{\Omega}(\ell_{ij})$: the link with the maximum $\overline{\Omega}(\ell_{ij})$; 
	\item[(5)] Undirected handling: first select the undirected link with the maximum $\Omega(\ell_{ij})$ by regarding the network is undirected, and then randomly determine the link direction.
\end{itemize}
Figure \ref{fig:performance_dir} shows that the method ``maximize $\Re(\mu)$" presents the best performance for a small number of adding links $K$ although the submodularity of the function $\Updelta\Re(\mu)$ in directed networks is not guaranteed. 
However, the generalized algebraic connectivity $\Re(\mu)$ via the method ``maximize $\Re(\mu)$" tends to increase monotonically but very slowly when the number of adding links is large, which implies that the original greedy method of ``maximize $\Re(\mu)$" could lead to the local optima in directed networks.
The method `` maximize $\underline{\Omega}(\ell_{ij})$ " is the best heuristic method to approach the method ``maximize $\Re(\mu)$" for a small fraction of adding links, and performs better for a large fraction of adding links.

We observe that the generalized algebraic connectivity $\Updelta\Re(\mu)$ via the method ``maximize $\underline{\Omega}(\ell_{ij})$" as a function of the number of adding links $K$ does not smoothly monotonically increase, which may instead help the solution to avoid local optimum possibly due to the inheritance of stochastic optimization (e.g., the principle of ``simulated annealing" algorithm and extremal optimization \cite{zhidong2013extremal}). 
Also, the method ``maximize $\overline{\Omega}(\ell_{ij})$ " can also avoid the local optimum and leads to a good solution, while the performance could degrade much in very sparse networks (e.g., Berlin traffic network with the average degree $E[D]\approx2.29$).

Figure \ref{fig:performance_dir} also shows the difference of performance between the method ``maximize $\underline{\Omega}(\ell_{ij})$ " and the same link with inverse direction, which demonstrates that the link direction is influential to the generalized algebraic connectivity $\Re(\mu)$. Furthermore, the unsatisfied performance of the method ``undirected handling'' hints the difference of the optimization between directed and undirected networks, as well as the fact that the strategy for undirected networks may be infeasible for directed networks.

\section{Conclusion}
The (generalized) algebraic connectivity indicates the lower-bound of the exponential convergence rate for consensus processes on networks. This paper investigated the strategy to maximize the (generalized) algebraic connectivity by adding links in both directed and undirected networks. 
The approximate submodularity of the derived perturbation of the algebraic connectivity validates that the greedy strategy could guarantee a constant performance in some undirected networks. For directed networks, we proposed the heuristic metrics for selecting the best directed links in the greedy algorithm based on the bounds of the generalized algebraic connectivity perturbation.
Numerical tests show that the original greedy strategy based on the exact algebraic connectivity performs best for a small number of adding links, while the greedy strategy based on the proposed metrics could outperform others when the number of adding links is large as well as requires less computational cost. 
This paper conveys the insight into the difference on the optimization in directed and undirected networks.
Beyond this paper, the similarities and differences of the behavior of dynamic processes in directed and undirected network merit further study, which is also of practical significance to the community detection for dynamics \cite{malliaros2013clustering}.

\appendices


\ifCLASSOPTIONcompsoc

\ifCLASSOPTIONcaptionsoff
  \newpage
\fi

\bibliographystyle{IEEEtran}
\bibliography{bibl}

\begin{thebibliography}{10}
\providecommand{\url}[1]{#1}
\csname url@samestyle\endcsname
\providecommand{\newblock}{\relax}
\providecommand{\bibinfo}[2]{#2}
\providecommand{\BIBentrySTDinterwordspacing}{\spaceskip=0pt\relax}
\providecommand{\BIBentryALTinterwordstretchfactor}{4}
\providecommand{\BIBentryALTinterwordspacing}{\spaceskip=\fontdimen2\font plus
\BIBentryALTinterwordstretchfactor\fontdimen3\font minus
  \fontdimen4\font\relax}
\providecommand{\BIBforeignlanguage}[2]{{%
\expandafter\ifx\csname l@#1\endcsname\relax
\typeout{** WARNING: IEEEtran.bst: No hyphenation pattern has been}%
\typeout{** loaded for the language `#1'. Using the pattern for}%
\typeout{** the default language instead.}%
\else
\language=\csname l@#1\endcsname
\fi
#2}}
\providecommand{\BIBdecl}{\relax}
\BIBdecl

\bibitem{van2014performance}
P.~Van~Mieghem, \emph{Performance analysis of complex networks and
  systems}.\hskip 1em plus 0.5em minus 0.4em\relax Cambridge University Press,
  2014.

\bibitem{simonsen2004diffusion}
I.~Simonsen, K.~A. Eriksen, S.~Maslov, and K.~Sneppen, ``Diffusion on complex
  networks: {A} way to probe their large-scale topological structures,''
  \emph{Physica A: Statistical Mechanics and its Applications}, vol. 336, no.
  1-2, pp. 163--173, 2004.

\bibitem{luyben1989process}
W.~L. Luyben, \emph{Process modeling, simulation and control for chemical
  engineers}.\hskip 1em plus 0.5em minus 0.4em\relax McGraw-Hill Higher
  Education, 1989.

\bibitem{van2017pseudoinverse}
P.~Van~Mieghem, K.~Devriendt, and H.~Cetinay, ``Pseudoinverse of the
  {L}aplacian and best spreader node in a network,'' \emph{Physical Review E},
  vol.~96, no.~3, p. 032311, 2017.

\bibitem{olfati2007consensus}
R.~Olfati-Saber, J.~A. Fax, and R.~M. Murray, ``Consensus and cooperation in
  networked multi-agent systems,'' \emph{Proceedings of the IEEE}, vol.~95,
  no.~1, pp. 215--233, 2007.

\bibitem{li2010consensus}
Z.~Li, Z.~Duan, G.~Chen, and L.~Huang, ``Consensus of multiagent systems and
  synchronization of complex networks: {A} unified viewpoint,'' \emph{IEEE
  Transactions on Circuits and Systems I: Regular Papers}, vol.~57, no.~1, pp.
  213--224, 2010.

\bibitem{ghosh2006growing}
A.~Ghosh and S.~Boyd, ``Growing well-connected graphs,'' in \emph{2006 45th
  IEEE Conference on Decision and Control}.\hskip 1em plus 0.5em minus
  0.4em\relax IEEE, 2006, pp. 6605--6611.

\bibitem{kim2006maximizing}
Y.~Kim and M.~Mesbahi, ``On maximizing the second smallest eigenvalue of a
  state-dependent graph {L}aplacian,'' \emph{IEEE Transactions on Automatic
  Control}, vol.~51, no.~1, pp. 116--120, 2006.

\bibitem{boyd2004fastest}
S.~Boyd, P.~Diaconis, and L.~Xiao, ``Fastest mixing {M}arkov chain on a
  graph,'' \emph{SIAM Review}, vol.~46, no.~4, pp. 667--689, 2004.

\bibitem{sydney2013optimizing}
A.~Sydney, C.~Scoglio, and D.~Gruenbacher, ``Optimizing algebraic connectivity
  by edge rewiring,'' \emph{Applied Mathematics and Computation}, vol. 219,
  no.~10, pp. 5465--5479, 2013.

\bibitem{ogiwara2017maximizing}
K.~Ogiwara, T.~Fukami, and N.~Takahashi, ``Maximizing algebraic connectivity in
  the space of graphs with a fixed number of vertices and edges,'' \emph{IEEE
  Transactions on Control of Network Systems}, vol.~4, no.~2, pp. 359--368,
  2017.

\bibitem{alenazi2013network}
M.~J. Alenazi, E.~K. Cetinkaya, and J.~P. Sterbenz, ``Network design and
  optimisation based on cost and algebraic connectivity,'' in \emph{2013 5th
  International Congress on Ultra Modern Telecommunications and Control Systems
  and Workshops (ICUMT)}.\hskip 1em plus 0.5em minus 0.4em\relax IEEE, 2013,
  pp. 193--200.

\bibitem{li2018maximizing}
G.~Li, Z.~F. Hao, H.~Huang, and H.~Wei, ``Maximizing algebraic connectivity via
  minimum degree and maximum distance,'' \emph{IEEE Access}, vol.~6, pp.
  41\,249--41\,255, 2018.

\bibitem{li2010cooperative}
C.~Li, Z.~Qu, A.~Das, and F.~Lewis, ``Cooperative control with improvable
  network connectivity,'' in \emph{American Control Conference (ACC),
  2010}.\hskip 1em plus 0.5em minus 0.4em\relax IEEE, 2010, pp. 87--92.

\bibitem{asadi2017expected}
M.~M. Asadi, M.~Khosravi, A.~G. Aghdam, and S.~Blouin, ``Expected convergence
  rate to consensus in asymmetric networks: {A}nalysis and distributed
  estimation,'' \emph{IEEE Transactions on Systems, Man, and Cybernetics:
  Systems}, 2017.

\bibitem{van2010graphspectra}
P.~Van~Mieghem, \emph{Graph spectra for complex networks}.\hskip 1em plus 0.5em
  minus 0.4em\relax Cambridge University Press, 2010.

\bibitem{stewart1990matrix}
G.~W. Stewart, ``Matrix perturbation theory,'' 1990.

\bibitem{li2013distributed}
C.~Li and Z.~Qu, ``Distributed estimation of algebraic connectivity of directed
  networks,'' \emph{Systems \& Control Letters}, vol.~62, no.~6, pp. 517--524,
  2013.

\bibitem{quarteroni2008matematica}
A.~Quarteroni, R.~Sacco, and F.~Saleri, \emph{Matematica numerica}.\hskip 1em
  plus 0.5em minus 0.4em\relax Springer Science \& Business Media, 2008.

\bibitem{milanese2010approximating}
A.~Milanese, J.~Sun, and T.~Nishikawa, ``Approximating spectral impact of
  structural perturbations in large networks,'' \emph{Physical Review E},
  vol.~81, no.~4, p. 046112, 2010.

\bibitem{lovasz1983submodular}
L.~Lov{\'a}sz, ``Submodular functions and convexity,'' in \emph{Mathematical
  Programming: The State of the Art}.\hskip 1em plus 0.5em minus 0.4em\relax
  Springer, 1983, pp. 235--257.

\bibitem{feldman2017greed}
M.~Feldman, C.~Harshaw, and A.~Karbasi, ``Greed is good: {N}ear-optimal
  submodular maximization via greedy optimization,'' \emph{arXiv preprint
  arXiv:1704.01652}, 2017.

\bibitem{nemhauser1978analysis}
G.~L. Nemhauser, L.~A. Wolsey, and M.~L. Fisher, ``An analysis of
  approximations for maximizing submodular set functions—{I},''
  \emph{Mathematical programming}, vol.~14, no.~1, pp. 265--294, 1978.

\bibitem{liu2018performance}
Y.~Liu, Z.~Zhang, E.~K. Chong, and A.~Pezeshki, ``Performance bounds with
  curvature for batched greedy optimization,'' \emph{Journal of Optimization
  Theory and Applications}, vol. 177, no.~2, pp. 535--562, 2018.

\bibitem{merris1994laplacian}
R.~Merris, ``Laplacian matrices of graphs: {A} survey,'' \emph{Linear Algebra
  and its Applications}, vol. 197, pp. 143--176, 1994.

\bibitem{gupta2018approximate}
G.~Gupta, S.~Pequito, and P.~Bogdan, ``Approximate submodular functions and
  performance guarantees,'' \emph{arXiv preprint arXiv:1806.06323}, 2018.

\bibitem{horel2016maximization}
T.~Horel and Y.~Singer, ``Maximization of approximately submodular functions,''
  in \emph{Advances in Neural Information Processing Systems}, 2016, pp.
  3045--3053.

\bibitem{van2010graph}
P.~Van~Mieghem, \emph{Graph spectra for complex networks}.\hskip 1em plus 0.5em
  minus 0.4em\relax Cambridge University Press, 2010.

\bibitem{kato2013perturbation}
T.~Kato, \emph{Perturbation theory for linear operators}.\hskip 1em plus 0.5em
  minus 0.4em\relax Springer Science \& Business Media, 2013, vol. 132.

\bibitem{ruhe1994rational}
A.~Ruhe, ``Rational {K}rylov algorithms for nonsymmetric eigenvalue problems,''
  in \emph{Recent Advances in Iterative Methods}.\hskip 1em plus 0.5em minus
  0.4em\relax Springer, 1994, pp. 149--164.

\bibitem{bauer1960norms}
F.~L. Bauer and C.~T. Fike, ``Norms and exclusion theorems,'' \emph{Numerische
  Mathematik}, vol.~2, no.~1, pp. 137--141, 1960.

\bibitem{belsley2005regression}
D.~A. Belsley, E.~Kuh, and R.~E. Welsch, \emph{Regression diagnostics:
  {I}dentifying influential data and sources of collinearity}.\hskip 1em plus
  0.5em minus 0.4em\relax John Wiley \& Sons, 2005, vol. 571.

\bibitem{zachary1977information}
W.~W. Zachary, ``An information flow model for conflict and fission in small
  groups,'' \emph{Journal of Anthropological Research}, vol.~33, no.~4, pp.
  452--473, 1977.

\bibitem{knuth1993stanford}
D.~E. Knuth, \emph{The {S}tanford {G}raph{B}ase: {A} platform for combinatorial
  computing}.\hskip 1em plus 0.5em minus 0.4em\relax Addison-Wesley Reading,
  1993, vol.~37.

\bibitem{newman2006finding}
M.~E.~J. Newman, ``Finding community structure in networks using the
  eigenvectors of matrices,'' \emph{Physical Review E}, vol.~74, no.~3, p.
  036104, 2006.

\bibitem{coleman1964introduction}
J.~S. Coleman, \emph{Introduction to mathematical sociology.}\hskip 1em plus
  0.5em minus 0.4em\relax New York: Free Press, 1964.

\bibitem{jahn2005system}
O.~Jahn, R.~H. M{\"o}hring, A.~S. Schulz, and N.~E. Stier-Moses,
  ``System-optimal routing of traffic flows with user constraints in networks
  with congestion,'' \emph{Operations Research}, vol.~53, no.~4, pp. 600--616,
  2005.

\bibitem{white1986structure}
J.~G. White, E.~Southgate, J.~N. Thomson, and S.~Brenner, ``The structure of
  the nervous system of the nematode {C}aenorhabditis elegans,'' \emph{Philos
  Trans R Soc Lond B Biol Sci}, vol. 314, no. 1165, pp. 1--340, 1986.

\bibitem{zhidong2013extremal}
H.~Zhidong and H.~Wenjun, ``Extremal optimization approach to joint routing and
  scheduling for industrial wireless networks,'' in \emph{2013 IEEE 10th
  International Conference on Mobile Ad-Hoc and Sensor Systems}.\hskip 1em plus
  0.5em minus 0.4em\relax IEEE, 2013, pp. 427--428.

\bibitem{malliaros2013clustering}
F.~D. Malliaros and M.~Vazirgiannis, ``Clustering and community detection in
  directed networks: {A} survey,'' \emph{Physics Reports}, vol. 533, no.~4, pp.
  95--142, 2013.

\end{thebibliography}

\end{document}